**Exact finite element formulation of quasi-3D beam based on analytical internal force fields for accurate static analysis of FG sandwich beams**

Wenxiong Li*, Suiyin Chen

*College of Water Conservancy and Civil Engineering, South China Agricultural University, Guangzhou 510642, China*

*Corresponding author. Email: leewenxiong@scau.edu.cn (W. Li)

**ABSTRACT**

This paper presents a novel exact finite element formulation of quasi-3D beam for high-fidelity analysis of functionally graded sandwich beams. Unlike conventional displacement-based elements that rely on approximate interpolation functions or existing exact finite element methods requiring closed-form solutions of generalized displacements, the proposed approach constructs an exact beam element from analytical expressions of internal forces. The method innovatively integrates force-based beam element formulation with the exact finite element framework through a systematic implementation procedure. Firstly, analytical expressions for internal forces are derived by applying differential equilibrium equations, geometric relations, and constitutive equations. These expressions are then used to obtain integral forms of generalized displacements expressed in terms of internal force parameters. Subsequently, nodal generalized displacements are defined, and their relationship with internal force parameters is established via the consistency condition. Finally, exact shape functions are derived from interpolation definitions, and the exact stiffness matrix and equilibrium equations are formulated using the principle of virtual work. To further enhance solution accuracy, modified cross-sectional stiffness matrices accounting for equilibrium-based stress distributions are incorporated. Comprehensive numerical examples demonstrate that the proposed element achieves superior solution accuracy and computational efficiency compared to conventional methods, validating its potential for precise structural analysis of complex beam structures with material gradients.

## 1. Introduction

Beams are extensively employed in diverse engineering sectors, from aerospace to civil infrastructure, owing to their straightforward mechanical design and exceptional applicability. In practice, many critical structural members can be idealized as beams, which remain indispensable elements even in complex systems predominantly composed of thin-walled assemblies. In addition, many specialized components, such as thin-walled beam-columns, are analyzed using theories derived from fundamental beam principles. Against this backdrop, Functionally Graded (FG) sandwich beams have attracted considerable interest as a highly promising structural element. By strategically tailoring their material composition and microstructure through the thickness, FG sandwich beams achieve significantly enhanced performance under complex loading and environmental conditions compared to their conventional counterparts. Consequently, with the growing adoption of FG sandwich beams in modern engineering systems, the development of highly accurate and computationally efficient models for their analysis has emerged as a key research focus.

To accurately predict the mechanical behavior of FG beams, a variety of beam theories have been developed, ranging from the Classical Beam Theory (CBT) [1-4] to shear deformation formulations such as the First-order Shear

Deformation Theory (FSDT) [5-13], Higher-order Shear Deformation Theory (HSDT) [14-27], and quasi-3D theories that incorporate both shear and transverse stretching deformation effects [28-36]. CBT neglects transverse shear deformation, rendering it unsuitable for thick beams where it tends to overestimate stiffness and underestimate deflection. FSDT accounts for shear deformation through a shear correction factor but relies on the assumption that the plane section remains plane, which becomes problematic for FG beams due to material heterogeneity. Moreover, FSDT fails to capture higher-order displacement variations across the thickness. HSDT employs higher-order displacement functions to model nonlinear displacement fields, eliminating the need for a shear correction factor and satisfying zero transverse shear stress at free surfaces. This advantage makes HSDT particularly effective for FG structures. In the early 1980s, Levinson [37] and Bickford [38] conducted pioneering work on HSDT, which formed the precursor of higher-order shear beam models. Subsequently, building upon their foundation, Reddy [39] established a more refined third-order shear deformation model. In recent years, HSDT has been increasingly applied to the analysis of FG sandwich beams, including the linear and nonlinear transient response analysis under moving distributed masses [40], the nonlinear bending behavior analysis based on different micromechanical models [41], and the post-buckling behavior analysis under high-temperature loads [42]. However, HSDT still neglects transverse normal stress and thickness-direction deformation, leading to inaccuracies in predicting principal stress distributions. To address these limitations, quasi-3D beam theories have been developed to account for both higher-order shear and transverse stretching deformation effects. Recently, the quasi-3D theories have been successfully applied to analyze the static, dynamic, nonlinear vibration, and buckling behaviors of various FG sandwich beams [28, 34, 35, 43]. Quasi-3D beam theory-based models have gained increasing interest because they not only capture essential features of stress and strain fields but also maintain reasonable computational efficiency. This dual capability enables more accurate predictions of principal stress distributions and deformation behaviors, making them particularly valuable for advanced structural analysis of FG materials.

Similar to conventional finite element-based beam models, the numerical implementation of quasi-3D beam theories typically relies on simplified interpolation functions to construct beam elements, leading to unavoidable discretization errors. These errors necessitate mesh refinement to achieve convergence, particularly in regions with high stress gradients or material transitions. For FG sandwich beams, the required level of mesh refinement depends strongly on the material distribution, making it difficult to determine an optimal mesh density systematically. This deficiency not only increases computational cost but also undermines the robustness and repeatability of numerical results, limiting the practical utility of quasi-3D beam models in engineering analysis. Consequently, the creation of a discretization-error-free quasi-3D beam element represents a breakthrough necessity, advancing theoretical frameworks while enabling reliable computational predictions for complex FG structures.

For displacement-based finite element models, the main approach to improving element accuracy is to construct a displacement interpolation function matrix that reflects the true response characteristics as much as possible. Advances in the exact finite element method have demonstrated its potential to achieve high precision in structural analysis [2, 19, 44-51]. This approach typically constructs finite elements using interpolation functions derived from the closed-form solutions of the differential equilibrium equations governing generalized displacement fields. Although conceptually attractive, this method faces a major limitation: the need to solve higher-order differential equilibrium equations governing displacement fields, which are often analytically intractable and computationally prohibitive. An alternative and increasingly explored strategy involves constructing finite elements based on closed-form solutions of internal force fields rather than displacement fields [10, 20, 52-56]. Since the differential equilibrium equations expressed in terms of internal forces generally exhibit a lower order than those for displacements, obtaining their analytical solutions is typically more feasible. This approach allows for the construction of beam elements with high precision by directly enforcing equilibrium conditions at the element boundaries and the compatibility condition within the element. Li et al. [20] successfully established a highly accurate finite element model for higher-order shear deformation beams using this method. However, this method requires the introduction of internal force parameters that describe internal force fields, resulting in more Degrees of Freedom (DoFs) for each element. Moreover, the other critical drawback is that the resulting stiffness matrix, while accurate, is often non-symmetric. This lack of symmetry complicates the numerical solution process and may reduce computational efficiency. Therefore, the establishment of an exact finite element method based on closed-form solutions of internal force fields and featuring a symmetric stiffness matrix, along with achieving high-precision solutions for quasi-3D beams through this methodology, represents a critical research frontier requiring further investigation. To the best of the authors' knowledge, integrating the exact finite element construction approach based on displacement formulations with the element development theory grounded in analytical internal force fields presents a breakthrough direction for novel exact finite element model development. This approach constructs exact quasi-3D beam elements by using low-order differential equilibrium equations governing internal force distributions and the interpolation function matrix's fundamental properties. It preserves stiffness matrix symmetry and cuts costs of solving higher-order equations, advancing high-performance finite element analysis without discretization errors.

For FG sandwich beams, accurate stress results serve as critical foundations for subsequent design and analysis procedures. However, most of the existing quasi-3D beam models fail to ensure the accuracy of stress predictions, particularly for transverse shear stress and transverse normal stress, which limits their practical application. This limitation is primarily due to the fact that the derivation of stress components in most quasi-3D models is based on kinematic assumptions and constitutive equations, without rigorous enforcement of the equilibrium differential equations. As a result, the predicted stress distributions may deviate significantly from the true solution, especially in FG sandwich beams where abrupt material property changes occur at layer interfaces. In response to this challenge, recent research has focused on developing enhanced beam theories and finite element models capable of accurately predicting stress distributions in FG sandwich beams. For transverse shear stress in higher-order beam models or quasi-3D beam models, various approaches have been proposed, including the global-local shear deformation theory [57, 58], mixed higher-order shear beam elements [21, 23, 30], and rational higher-order displacement functions based on equilibrium constraints [18, 22]. These methods aim to improve the accuracy of transverse shear stress predictions by incorporating additional kinematic or equilibrium conditions. As for transverse normal stress in quasi-3D beam models, several computational strategies have been explored. Salami et al. [59] developed a higher-order sandwich panel theory to improve stress predictions in moderately thick plates, Nguyen et al. [60] introduced a quasi-3D quadrilateral element for FG porous sandwich plates, and Van Vihn [61] formulated a hybrid quasi-3D theory combining polynomial and trigonometric functions for bi-functionally graded plates. These methods enhance stress accuracy either by increasing the number of DoFs or through post-processing corrections. However, they are not without drawbacks: increasing DoFs leads to higher computational costs, while the introduction of post-processing techniques often fails to ensure the accuracy of displacement solutions. Consequently, the development of a comprehensive quasi-3D beam model capable of simultaneously providing accurate solutions for both transverse shear stress, transverse normal stress, and displacement fields remains an open challenge in structural mechanics. A key approach involves directly formulating stress expressions from the fundamental differential equilibrium equations governing stress interrelationships, and integrating these with exact finite element construction methods.

This paper develops a quasi-3D beam exact finite element method based on the closed-form solutions of internal force fields, aiming at achieving high-fidelity analysis of FG sandwich beams. Unlike the conventional displacement-based element models that rely on approximate displacement interpolation functions or the exact finite element methods that require closed-form solutions of generalized displacements, the proposed approach constructs an exact beam element based on analytical expressions of the internal force fields. The proposed method innovatively integrates the force-based beam element formulation with the exact finite element framework, and the main implementation steps of this method are outlined as: (a) The analytical expressions for the internal forces are derived by systematically applying the differential equilibrium equations, geometric relations and cross-sectional constitutive equations of the quasi-3D beam; (b) The integral forms of generalized displacement fields are obtained based on the internal force expressions and the beam's geometric relations, expressed in terms of the internal force parameters and the generalized displacements at the starting node of the element; (c) The generalized displacements for nodes are defined, and the relationship between the internal force parameters and the nodal generalized displacements is formulated by utilizing the consistency condition; (d) The exact shape function matrices are derived from the interpolation definition, and the exact stiffness matrix and equilibrium equations of the quasi-3D beam element are obtained according to the principle of virtual work. In addition, to enhance the solution accuracy for FG sandwich beams, the proposed element incorporates the modified cross-sectional stiffness matrices, which account for the influence of the equilibrium-based stress distributions. The rest of this paper is organized as follows. Section 2 presents the quasi-3D theory for beam analysis with FG material. Section 3 details the formulation of the force-based exact finite element model, deriving the exact stiffness matrix from the internal force field expressions. Section 4 describes the method implemented to enhance the accuracy of stress solutions. The accuracy and efficiency of the proposed beam element model are validated through numerical examples in Section 5. Finally, Section 6 summarizes the conclusions of this study.

## 2. Quasi-3D theory for beam analysis

For the FG sandwich beam with length $L$ and cross-section $b \times h$, as shown in **Fig. 1**, the displacement fields can be expressed by generalized displacements $u(x)$, $w(x)$, $\theta(x)$ and $\phi(x)$ as follows based on the quasi-3D beam theory

$$u_x(x,y) = u(x) - y\frac{\mathrm{d}w(x)}{\mathrm{d}x} + f(y)\theta(x) \tag{1a}$$

$$u_y(x,y) = w(x) + g(y)\phi(x) \tag{1b}$$

where $u_x(x,y)$ and $u_y(x,y)$ are the displacement components of the point ($x$, $y$) along with $x$ and $y$ directions, respectively, $u(x)$ and $w(x)$ represent the axial displacement and transverse displacement of the beam axis, respectively, $\theta(x)$ refers to the magnitude of higher-order axial displacement, $\phi(x)$ is the measurement of the additional transverse displacement from transverse stretching deformation. In Eq. (1), $f(y)$ and $g(y)$ are the higher-order functions describing the distribution of out-of-plane axial displacement through the thickness and the distribution of the additional transverse displacement, respectively. They can be determined by different mathematical functions, including Third-order polynomial function (T), Sinusoidal function (S), and Exponential function (E), as shown in

**Table 1**. Owing to their capacity to inherently satisfy the stress free boundary conditions at the top and bottom surfaces, these three types of functions have garnered significant attention and are widely employed in related studies [28, 35, 41, 43].

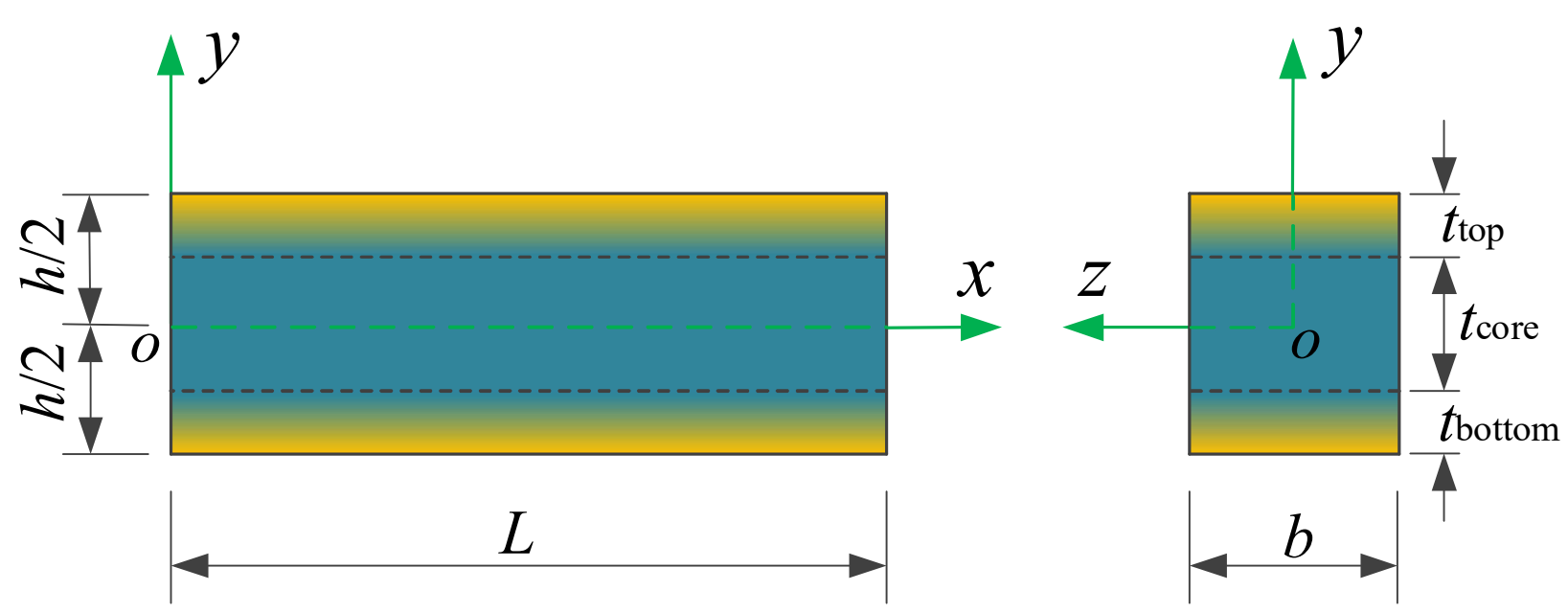

**Fig. 1.** FG sandwich beam with length $L$ and cross-section $b \times h$.

**Table 1** Functions for determining higher-order components.

| Function form | $f(y)$ | $g(y)$ | $\mathrm{d}g(y)/\mathrm{d}y$ |
|---|---|---|---|
| Third-order polynomial (T) | $y-\dfrac{4y^3}{3h^2}$ | $1-\dfrac{4y^2}{h^2}$ | $-\dfrac{8y}{h^2}$ |
| Sinusoidal (S) | $\dfrac{h}{\pi}\sin\left(\dfrac{\pi y}{h}\right)$ | $\cos\left(\dfrac{\pi y}{h}\right)$ | $-\dfrac{\pi}{h}\sin\left(\dfrac{\pi y}{h}\right)$ |
| Exponential (E) | $ye^{\left(-2y^2/h^2\right)}$ | $\left(1-4y^2/h^2\right)e^{\left(-2y^2/h^2\right)}$ | $\left(16y^3/h^4-12y/h^2\right)e^{\left(-2y^2/h^2\right)}$ |

From Eq. (1), the geometric relations defining the expressions of strain components, including axial normal strain, transverse normal strain and transverse shear strain, can be derived as

$$\varepsilon_x(x,y)=\frac{\partial u_x(x,y)}{\partial x}=\mathbf{B}_x(y)\boldsymbol{\varepsilon}_\mathrm{n}(x) \tag{2a}$$

$$\varepsilon_y(x,y)=\frac{\partial u_y(x,y)}{\partial y}=\mathbf{B}_y(y)\boldsymbol{\varepsilon}_\mathrm{n}(x) \tag{2b}$$

$$\gamma_{xy}(x,y)=\frac{\partial u_x(x,y)}{\partial y}+\frac{\partial u_y(x,y)}{\partial x}=g(y)\left[\theta(x)+\phi_{,x}(x)\right] \tag{2c}$$

where $\boldsymbol{\varepsilon}_\mathrm{n}(x)$ represents the generalized strain vector related to the normal strains, which is expressed as

$$\boldsymbol{\varepsilon}_\mathrm{n}(x)=\left\{u_{,x}(x) \quad -w_{,xx}(x) \quad \theta_{,x}(x) \quad \phi(x)\right\}^\mathrm{T}, \tag{3}$$

and $\mathbf{B}_x(y)$ and $\mathbf{B}_y(y)$ are expressed as

$$\mathbf{B}_x(y)=\left\{1 \quad y \quad f(y) \quad 0\right\} \tag{4a}$$

$$\mathbf{B}_y(y)=\left\{0 \quad 0 \quad 0 \quad g_{,y}(y)\right\} \tag{4b}$$

Note that $(\cdot)_{,x}=\mathrm{d}(\cdot)/\mathrm{d}x$ and $(\cdot)_{,xx}=\mathrm{d}^2(\cdot)/\mathrm{d}x^2$ are applied to indicate the first-order and second-order differentiations with respect to coordinate $x$, respectively, $(\cdot)_{,y}=\mathrm{d}(\cdot)/\mathrm{d}y$ is applied to indicate the first-order differentiation with respect to coordinate $y$, and the superscript 'T' denotes the transpose of the vector or matrix, as shown in Eq. (3).

For FG beams where the material properties vary through the thickness, the stress-strain relationship for the quasi-3D model can be expressed as

$$\begin{bmatrix}\sigma_x(x,y)\\ \sigma_y(x,y)\\ \tau_{xy}(x,y)\end{bmatrix}=\begin{bmatrix}\bar{Q}_{11}(y) & \bar{Q}_{12}(y) & 0\\ \bar{Q}_{12}(y) & \bar{Q}_{11}(y) & 0\\ 0 & 0 & \bar{Q}_{33}(y)\end{bmatrix}\begin{bmatrix}\varepsilon_x(x,y)\\ \varepsilon_y(x,y)\\ \gamma_{xy}(x,y)\end{bmatrix} \tag{5}$$

where $\bar{Q}_{11}(y)=\dfrac{E(y)}{1-v^2(y)}$, $\bar{Q}_{12}(y)=\dfrac{v(y)E(y)}{1-v^2(y)}$ and $\bar{Q}_{33}(y)=\dfrac{E(y)}{2\left[1+v(y)\right]}$ with $E(y)$ and $v(y)$ being the Young's modulus and Poisson's ratio, respectively.

Based on the strain expressions in Eq. (2), the variation of strain energy can be expressed as

$$\begin{aligned}\delta U=&\int_0^L\int_A\sigma_x(x,y)\left[\delta u_{,x}(x)-y\delta w_{,xx}(x)+f(y)\delta\theta_{,x}(x)\right]\mathrm{d}A\mathrm{d}x+\int_0^L\int_A\sigma_y(x,y)g_{,y}(y)\delta\phi(x)\mathrm{d}A\mathrm{d}x\\ &+\int_0^L\int_A\tau_{xy}(x,y)g(y)\left[\delta\theta(x)+\delta\phi_{,x}(x)\right]\mathrm{d}A\mathrm{d}x\end{aligned} \tag{6}$$

where $L$ represents the length of the beam and $A$ denotes the domain of the beam's cross-section. By integrating across the beam's cross-section, the expression of variation of strain energy can be expressed as

$$\delta U=\int_0^L\left[N(x)\delta u_{,x}(x)-M_w(x)\delta w_{,xx}(x)+M_\theta(x)\delta\theta_{,x}(x)+Q(x)\delta\theta(x)+Q(x)\delta\phi_{,x}(x)+R(x)\delta\phi(x)\right]\mathrm{d}x \tag{7}$$

where $N(x)$ is the axial force of the beam; $M_w(x)$ and $M_\theta(x)$ are the bending moments corresponding to the first-order and higher-order curvatures, respectively; $Q(x)$ is the shear force corresponding to the higher-order shear deformation and $R(x)$ is the stress resultant corresponding to the transverse stretching deformation. These internal force components are defined as

$$N(x)=\int_A\sigma_x(x,y)\mathrm{d}A \tag{8a}$$

$$M_w(x)=\int_A y\sigma_x(x,y)\mathrm{d}A \tag{8b}$$

$$M_\theta(x)=\int_A f(y)\sigma_x(x,y)\mathrm{d}A \tag{8c}$$

$$Q(x)=\int_A g(y)\tau_{xy}(x,y)\mathrm{d}A \tag{8d}$$

$$R(x)=\int_A g_{,y}(y)\sigma_y(x,y)\mathrm{d}A \tag{8e}$$

For bending analysis with distributed load $q(x)$ on the upper surface ($y=h/2$) and the nodal forces at the two ending nodes, as shown in **Fig. 2**, the variation of potential energy can be expressed as

$$\delta V=-\int_0^L q(x)\delta w(x)\mathrm{d}x-P_{ax}\delta u(0)-P_{ay}\delta w(0)-M_a\delta w_{,x}(0)-P_{bx}\delta u(L)-P_{by}\delta w(L)-M_b\delta w_{,x}(L) \tag{9}$$

Then, the following equation can be established according to $\delta U+\delta V=0$

$$\begin{aligned}&\int_0^L\left[N(x)\delta u_{,x}(x)-M_w(x)\delta w_{,xx}(x)+M_\theta(x)\delta\theta_{,x}(x)+Q(x)\delta\theta(x)+Q(x)\delta\phi_{,x}(x)+R(x)\delta\phi(x)\right]\mathrm{d}x\\ &-\int_0^L q(x)\delta w(x)\mathrm{d}x-P_{ax}\delta u(0)-P_{ay}\delta w(0)-M_a\delta w_{,x}(0)-P_{bx}\delta u(L)-P_{by}\delta w(L)-M_b\delta w_{,x}(L)=0\end{aligned} \tag{10}$$

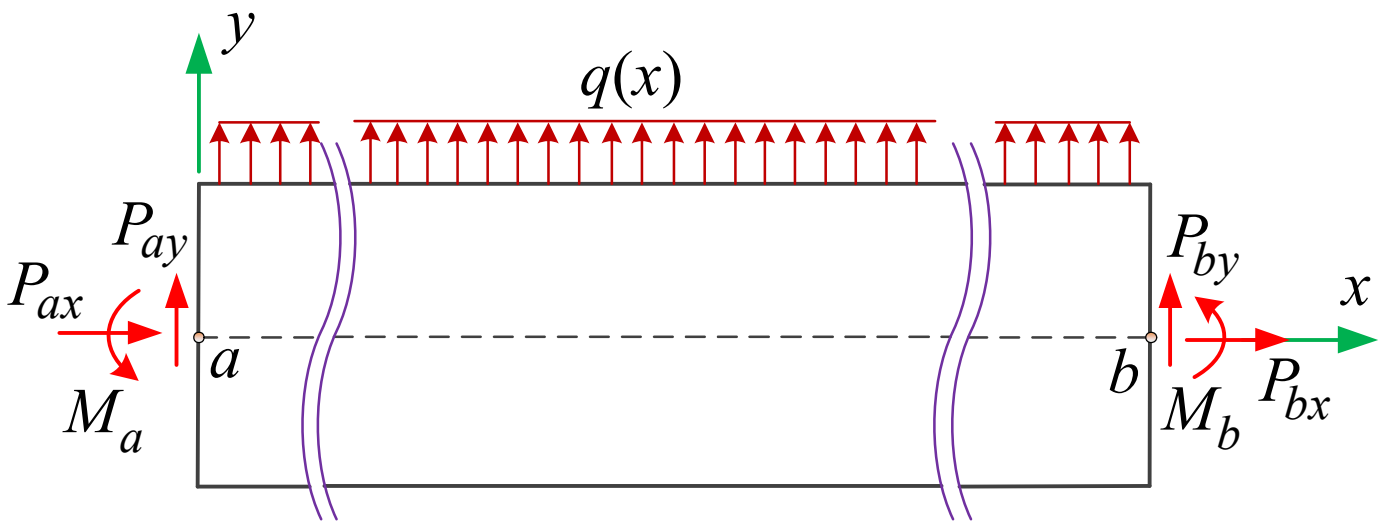

**Fig. 2.** Nodal loads and transverse distributed load at $y = h/2$.

Considering the arbitrariness of the virtual displacements, the following differential equilibrium equations governing internal force fields and the corresponding boundary conditions can be obtained from Eq. (10):

(1) Differential equilibrium equations

$$N_{,x}(x) = 0 \tag{11a}$$

$$M_{w,xx}(x) + q(x) = 0 \tag{11b}$$

$$M_{\theta,x}(x) - Q(x) = 0 \tag{11c}$$

$$Q_{,x}(x) - R(x) = 0 \tag{11d}$$

(2) Boundary conditions

at $x = 0$:

$$\delta u(0) = 0 \quad \text{or} \quad N(0) + P_{ax} = 0 \tag{12a}$$

$$\delta w(0) = 0 \quad \text{or} \quad M_{w,x}(0) + P_{ay} = 0 \tag{12b}$$

$$\delta w_{,x}(0) = 0 \quad \text{or} \quad M_{w}(0) + M_{a} = 0 \tag{12c}$$

$$\delta \theta(0) = 0 \quad \text{or} \quad M_{\theta}(0) = 0 \tag{12d}$$

$$\delta \phi(0) = 0 \quad \text{or} \quad Q(0) = 0 \tag{12e}$$

at $x = L$:

$$\delta u(L) = 0 \quad \text{or} \quad N(L) - P_{bx} = 0 \tag{13a}$$

$$\delta w(L) = 0 \quad \text{or} \quad M_{w,x}(L) - P_{by} = 0 \tag{13b}$$

$$\delta w_{,x}(L) = 0 \quad \text{or} \quad M_{w}(L) - M_{b} = 0 \tag{13c}$$

$$\delta \theta(L) = 0 \quad \text{or} \quad M_{\theta}(L) = 0 \tag{13d}$$

$$\delta \phi(L) = 0 \quad \text{or} \quad Q(L) = 0 \tag{13e}$$

By substituting Eqs. (5) and (2) into Eq. (8), the constitutive equations of the beam's cross-section can be obtained as

$$\boldsymbol{\sigma}_{\mathrm{n}}(x) = \mathbf{D}_{\mathrm{n}} \boldsymbol{\varepsilon}_{\mathrm{n}}(x) \tag{14}$$

$$Q(x) = D_{55}\left[\theta(x) + \phi_{,x}(x)\right] \tag{15}$$

where $\boldsymbol{\sigma}_{\mathrm{n}}(x)$ and $\mathbf{D}_{\mathrm{n}}$ are considered as the internal force vector and cross-section stiffness matrix related to normal stresses, respectively, and they are expressed as

$$\boldsymbol{\sigma}_{\mathrm{n}}(x)=\begin{Bmatrix} N(x) \\ M_w(x) \\ M_\theta(x) \\ R(x) \end{Bmatrix},\ \mathbf{D}_{\mathrm{n}}=\begin{bmatrix} D_{11} & D_{12} & D_{13} & D_{14} \\ D_{12} & D_{22} & D_{23} & D_{24} \\ D_{13} & D_{23} & D_{33} & D_{34} \\ D_{14} & D_{24} & D_{34} & D_{44} \end{bmatrix} \tag{16}$$

Specifically, the stiffness coefficients $D_{ij}\left(i,j=1,2,3,4\right)$ in $\mathbf{D}_{\mathrm{n}}$ and $D_{55}$ can be obtained by

$$\left[D_{11},D_{12},D_{22}\right]=\int_A \bar{Q}_{11}(y)\left[1,y,y^2\right]\mathrm{d}A \tag{17a}$$

$$\left[D_{13},D_{23},D_{33}\right]=\int_A \bar{Q}_{11}(y)f(y)\left[1,y,f(y)\right]\mathrm{d}A \tag{17b}$$

$$\left[D_{14},D_{24},D_{34}\right]=\int_A \bar{Q}_{12}(y)g_{,y}(y)\left[1,y,f(y)\right]\mathrm{d}A \tag{17c}$$

$$D_{44}=\int_A \bar{Q}_{11}(y)g_{,y}^2(y)\mathrm{d}A \tag{17d}$$

$$D_{55}=\int_A \bar{Q}_{33}(y)g^2(y)\mathrm{d}A \tag{17e}$$

By substituting Eqs. (14) and (15) into Eq. (11), the differential equilibrium equations governing generalized displacement fields can be obtained as

$$D_{11}u_{,xx}(x)-D_{12}w_{,xxx}(x)+D_{13}\theta_{,xx}(x)+D_{14}\phi_{,x}(x)=0 \tag{18a}$$

$$D_{12}u_{,xxx}(x)-D_{22}w_{,xxxx}(x)+D_{23}\theta_{,xxx}(x)+D_{24}\phi_{,xx}(x)+q(x)=0 \tag{18b}$$

$$D_{31}u_{,xx}(x)-D_{32}w_{,xxx}(x)+D_{33}\theta_{,xx}(x)-D_{55}\theta(x)+\left(D_{34}-D_{55}\right)\phi_{,x}(x)=0 \tag{18c}$$

$$D_{11}u_{,x}(x)-D_{12}w_{,xx}(x)+\left(D_{13}-D_{55}\right)\theta_{,x}(x)-D_{55}\phi_{,xx}(x)+D_{14}\phi(x)=0 \tag{18d}$$

Generally, the exact finite element method directly derives closed-form solutions of generalized displacements from Eq. (18) to construct exact shape function matrices, which are subsequently used to formulate exact stiffness matrix. However, the challenge of solving 12th-order differential equations in Eq. (18) motivates our alternative approach: leveraging the tractable lower-order system (Eq. (11)) to establish analytical solutions of internal forces, thereby providing a mathematically rigorous foundation for quasi-3D beam exact finite element formulation.

## 3. Formulation of force-based exact finite element

### *3.1. Closed-form solutions of internal forces*

For the sake of simplicity, the distributed load is expressed by the quadratic polynomial function as

$$q(x)=q_0+q_1x+q_2x^2 \tag{19}$$

with $q_0$, $q_1$ and $q_2$ representing the load parameters. According to the first and second equations in Eq. (11), the expressions of $N(x)$ and $M_w(x)$ can be obtained as

$$N(x)=c_0 \tag{20}$$

$$M_w(x)=c_1+c_2x-\frac{q_0x^2}{2}-\frac{q_1x^3}{6}-\frac{q_2x^4}{12} \tag{21}$$

where $c_0$, $c_1$ and $c_2$ are the internal force parameters to be determined.

Considering that the remaining two equations (the third and fourth equations) in Eq. (11) are not sufficient to determine the solutions of $M_\theta(x)$, $Q(x)$ and $R(x)$, an additional equation from the constitutive and geometric relations is required to achieve a solution.

Firstly, according to Eq. (14), the generalized strain vector can be expressed by the internal force vector as

$$\boldsymbol{\varepsilon}_{\mathrm{n}}(x)=\mathbf{F}_{\mathrm{n}}\boldsymbol{\sigma}_{\mathrm{n}}(x) \tag{22}$$

where $\mathbf{F}_{\mathrm{n}}$ is the cross-section flexibility matrix related to normal stresses, which is expressed as

$$\mathbf{F}_{\mathrm{n}}=\begin{bmatrix} f_{11} & f_{12} & f_{13} & f_{14} \\ f_{12} & f_{22} & f_{23} & f_{24} \\ f_{13} & f_{23} & f_{33} & f_{34} \\ f_{14} & f_{24} & f_{34} & f_{44} \end{bmatrix}=\mathbf{D}_{\mathrm{n}}^{-1} \tag{23}$$

Then, the following two equations can be obtained

$$\theta_{,x}(x)=f_{13}N(x)+f_{23}M_w(x)+f_{33}M_\theta(x)+f_{34}R(x) \tag{24}$$

$$\phi(x)=f_{14}N(x)+f_{24}M_w(x)+f_{34}M_\theta(x)+f_{44}R(x) \tag{25}$$

By substituting Eqs. (24) and (25) into Eq. (15), the following additional equation can be derived

$$f_{13}N(x)+f_{23}M_w(x)+f_{33}M_\theta(x)+f_{34}R(x)+f_{14}N_{,xx}(x)+f_{24}M_{w,xx}(x)+f_{34}M_{\theta,xx}(x)+f_{44}R_{,xx}(x)=D_{55}^{-1}Q_{,x}(x) \tag{26}$$

By using the third and fourth equations in Eq. (11a) and introducing the expressions of $N(x)$ and $M_w(x)$ (see Eqs. (20) and (21)), Eq. (26) can be further rewritten as

$$\begin{aligned} & M_{\theta,xxxx}(x)+b_1M_{\theta,xx}(x)+b_0M_\theta(x) \\ & =\frac{f_{23}}{12f_{44}}q_2x^4+\frac{f_{23}}{6f_{44}}q_1x^3+\left(\frac{f_{23}}{2f_{44}}q_0+\frac{f_{24}}{f_{44}}q_2\right)x^2+\left(\frac{f_{24}}{f_{44}}q_1-\frac{f_{23}}{f_{44}}c_2\right)x+\frac{f_{24}}{f_{44}}q_0-\frac{f_{13}}{f_{44}}c_0-\frac{f_{23}}{f_{44}}c_1 \end{aligned} \tag{27}$$

where

$$b_0=\frac{f_{33}}{f_{44}},\ b_1=\frac{2f_{34}-D_{55}^{-1}}{f_{44}} \tag{28}$$

Only a fourth-order linear ordinary differential equation needs to be solved to obtain the analytical solutions of the internal force fields. This approach significantly simplifies the process compared to directly solving for the displacement fields.

The eigen-equation corresponding to the homogeneous equation of Eq. (27) can be represented as

$$\lambda^4+b_1\lambda^2+b_0=0 \tag{29}$$

By solving the equation above, two pairs of conjugate complex roots with opposite real parts and identical imaginary parts are generally obtained, which are expressed as

$$\lambda_{1,2}=\alpha\pm\beta i \tag{30a}$$

$$\lambda_{3,4}=-\alpha\pm\beta i \tag{30b}$$

where $\alpha$ and $\beta$ represent the real part coefficient and imaginary part coefficient, respectively. Then, the general solution of the equation can be expressed as

$$\bar{M}_\theta(x)=c_3e^{\alpha x}\cos(\beta x)+c_4e^{\alpha x}\sin(\beta x)+c_5e^{-\alpha x}\cos(\beta x)+c_6e^{-\alpha x}\sin(\beta x) \tag{31}$$

where $c_3\sim c_6$ are the internal force parameters to be determined.

Further, the particular solution of Eq. (27) can be predefined as

$$M_\theta^*\left(x\right)=a_4x^4+a_3x^3+a_2x^2+a_1x+a_0 \tag{32}$$

The introduction of Eq. (32) into Eq. (27) can derive

$$\begin{aligned}&\left(b_0a_4-\frac{q_2f_{23}}{12f_{44}}\right)x^4+\left(b_0a_3-\frac{q_1f_{23}}{6f_{44}}\right)x^3+\left(b_0a_2+12b_1a_4-\frac{q_0f_{23}}{2f_{44}}-\frac{q_2f_{24}}{f_{44}}\right)x^2\\&+\left(6b_1a_3+b_0a_1-\frac{f_{24}}{f_{44}}q_1+\frac{f_{23}}{f_{44}}c_2\right)x+\left(24a_4+2b_1a_2+b_0a_0-\frac{f_{24}}{f_{44}}q_0+\frac{f_{13}}{f_{44}}c_0+\frac{f_{23}}{f_{44}}c_1\right)=0\end{aligned} \tag{33}$$

To ensure the constancy of Eq. (33), $a_0 \sim a_4$ are taken as

$$a_0=\frac{q_0f_{24}}{b_0f_{44}}-\frac{q_0b_1f_{23}}{b_0^2f_{44}}-\frac{f_{13}}{b_0f_{44}}c_0-\frac{f_{23}}{b_0f_{44}}c_1-\frac{2q_2f_{23}}{b_0^2f_{44}}-\frac{2q_2b_1f_{24}}{b_0^2f_{44}}+\frac{2q_2b_1^2f_{23}}{b_0^3f_{44}} \tag{34a}$$

$$a_1=\frac{q_1f_{24}}{b_0f_{44}}-\frac{q_1b_1f_{23}}{b_0^2f_{44}}-\frac{f_{23}}{b_0f_{44}}c_2 \tag{34b}$$

$$a_2=\frac{q_0f_{23}}{2b_0f_{44}}+\frac{q_2f_{24}}{b_0f_{44}}-\frac{q_2b_1f_{23}}{b_0^2f_{44}} \tag{34c}$$

$$a_3=\frac{q_1f_{23}}{6b_0f_{44}} \tag{34d}$$

$$a_4=\frac{q_2f_{23}}{12b_0f_{44}} \tag{34e}$$

Hence, the closed-form solution of $M_\theta\left(x\right)$ is

$$\begin{aligned}M_\theta\left(x\right)=&\,g_0c_0+g_1c_1+g_1xc_2+z_3\left(x\right)c_3+z_4\left(x\right)c_4+z_5\left(x\right)c_5+z_6\left(x\right)c_6\\&+\frac{h_2q_2x^4}{12}+\frac{h_2q_1x^3}{6}+\left(\frac{h_2q_0}{2}+\frac{k_2q_2}{2}\right)x^2+h_0q_1x+h_0q_0+k_0q_2\end{aligned} \tag{35}$$

where

$$z_3\left(x\right)=e^{\alpha x}\cos\left(\beta x\right) \tag{36a}$$

$$z_4\left(x\right)=e^{\alpha x}\sin\left(\beta x\right) \tag{36b}$$

$$z_5\left(x\right)=e^{-\alpha x}\cos\left(\beta x\right) \tag{36c}$$

$$z_6\left(x\right)=e^{-\alpha x}\sin\left(\beta x\right) \tag{36d}$$

and

$$g_0=-\frac{f_{13}}{b_0f_{44}},\ g_1=-\frac{f_{23}}{b_0f_{44}} \tag{37a}$$

$$h_0=\frac{f_{24}}{b_0f_{44}}-\frac{b_1f_{23}}{b_0^2f_{44}},\ h_2=\frac{f_{23}}{b_0f_{44}} \tag{37b}$$

$$k_0=-\frac{2f_{23}}{b_0^2f_{44}}-\frac{2b_1f_{24}}{b_0^2f_{44}}+\frac{2b_1^2f_{23}}{b_0^3f_{44}},\ k_2=\frac{2f_{24}}{b_0f_{44}}-\frac{2b_1f_{23}}{b_0^2f_{44}} \tag{37c}$$

By using Eq. (11) that $Q(x)=M_{\theta,x}(x)$ and $R(x)=Q_{,x}(x)$, the solutions of $Q(x)$ and $R(x)$ can also be derived. Then, the closed-form solutions of all five internal force components are characterized by seven internal force parameters ($c_0 \sim c_6$) and can be expressed as

$$\begin{Bmatrix} N(x) \\ M_w(x) \\ M_\theta(x) \\ R(x) \\ Q(x) \end{Bmatrix} = \begin{bmatrix} 1 & 0 & 0 & 0 & 0 & 0 & 0 \\ 0 & 1 & x & 0 & 0 & 0 & 0 \\ g_0 & g_1 & g_1 x & z_3(x) & z_4(x) & z_5(x) & z_6(x) \\ 0 & 0 & 0 & z_{3,xx}(x) & z_{4,xx}(x) & z_{5,xx}(x) & z_{6,xx}(x) \\ 0 & 0 & g_1 & z_{3,x}(x) & z_{4,x}(x) & z_{5,x}(x) & z_{6,x}(x) \end{bmatrix} \mathbf{c} + \begin{Bmatrix} 0 \\ -q_2 \\ h_2 q_2 \\ 0 \\ 0 \end{Bmatrix} \frac{x^4}{12} + \begin{Bmatrix} 0 \\ -q_1 \\ h_2 q_1 \\ 0 \\ 2h_2 q_2 \end{Bmatrix} \frac{x^3}{6} + \begin{Bmatrix} 0 \\ -q_0 \\ h_2 q_0 + k_2 q_2 \\ 2h_2 q_2 \\ h_2 q_1 \end{Bmatrix} \frac{x^2}{2} + \begin{Bmatrix} 0 \\ 0 \\ h_0 q_1 \\ h_2 q_1 \\ h_2 q_0 + k_2 q_2 \end{Bmatrix} x + \begin{Bmatrix} 0 \\ 0 \\ h_0 q_0 + k_0 q_2 \\ h_2 q_0 + k_2 q_2 \\ h_0 q_1 \end{Bmatrix} \tag{38}$$

where $\mathbf{c}$ represents the internal force parameter vector, which is expressed as

$$\mathbf{c} = \{c_0 \quad c_1 \quad c_2 \quad c_3 \quad c_4 \quad c_5 \quad c_6\}^{\mathrm{T}} \tag{39}$$

Note that the expressions of $z_{3,x}(x) \sim z_{6,x}(x)$ and $z_{3,xx}(x) \sim z_{6,xx}(x)$ can be explicitly derived according to the expressions of $z_3(x) \sim z_6(x)$.

*3.2. Expressions of generalized displacement fields*

Since the internal force expressions have been established, the constitutive laws of cross-section (Eqs. (22) and (15)) can be utilized to derive the generalized strains mathematically and subsequently the integral expressions for each generalized displacement component can be derived.

For the sake of convenience, the internal force vector related to normal stresses can be expressed by the internal force parameters as

$$\boldsymbol{\sigma}_{\mathrm{n}}(x) = \mathbf{N}_\sigma(x)\mathbf{c} + \mathbf{F}_\sigma(x) \tag{40}$$

where

$$\mathbf{N}_\sigma(x) = \begin{bmatrix} 1 & 0 & 0 & 0 & 0 & 0 & 0 \\ 0 & 1 & x & 0 & 0 & 0 & 0 \\ g_0 & g_1 & g_1 x & z_3(x) & z_4(x) & z_5(x) & z_6(x) \\ 0 & 0 & 0 & z_{3,xx}(x) & z_{4,xx}(x) & z_{5,xx}(x) & z_{6,xx}(x) \end{bmatrix} \tag{41}$$

$$\mathbf{F}_\sigma(x) = \begin{Bmatrix} 0 \\ -q_2 \\ h_2 q_2 \\ 0 \end{Bmatrix} \frac{x^4}{12} + \begin{Bmatrix} 0 \\ -q_1 \\ h_2 q_1 \\ 0 \end{Bmatrix} \frac{x^3}{6} + \begin{Bmatrix} 0 \\ -q_0 \\ h_2 q_0 + k_2 q_2 \\ 2h_2 q_2 \end{Bmatrix} \frac{x^2}{2} + \begin{Bmatrix} 0 \\ 0 \\ h_0 q_1 \\ h_2 q_1 \end{Bmatrix} x + \begin{Bmatrix} 0 \\ 0 \\ h_0 q_0 + k_1 q_2 \\ h_2 q_0 + k_2 q_2 \end{Bmatrix} \tag{42}$$

Meanwhile, the shear force is expressed as

$$Q(x) = \mathbf{N}_Q(x)\mathbf{c} + \mathbf{F}_Q(x) \tag{43}$$

where

$$\mathbf{N}_Q(x) = \begin{bmatrix} 0 & 0 & g_1 & z_{3,x}(x) & z_{4,x}(x) & z_{5,x}(x) & z_{6,x}(x) \end{bmatrix} \tag{44}$$

$$\mathbf{F}_Q(x) = \frac{2h_2q_2}{6}x^3 + \frac{h_2q_1}{2}x^2 + (h_2q_0 + k_2q_2)x + h_0q_1 \tag{45}$$

Then, the differentiations of generalized displacements can be obtained according to Eqs. (22) and (15) as

$$u_{,x}(x) = \mathbf{T}_u\mathbf{F}_\mathrm{n}\left[\mathbf{N}_\sigma(x)\mathbf{c} + \mathbf{F}_\sigma(x)\right] \tag{46}$$

$$w_{,xx}(x) = \mathbf{T}_w\mathbf{F}_\mathrm{n}\left[\mathbf{N}_\sigma(x)\mathbf{c} + \mathbf{F}_\sigma(x)\right] \tag{47}$$

$$\theta_{,x}(x) = \mathbf{T}_\theta\mathbf{F}_\mathrm{n}\left[\mathbf{N}_\sigma(x)\mathbf{c} + \mathbf{F}_\sigma(x)\right] \tag{48}$$

$$\phi_{,x}(x) = D_{55}^{-1}\left[\mathbf{N}_Q(x)\mathbf{c} + \mathbf{F}_Q(x)\right] - \theta(x) \tag{49}$$

where

$$\mathbf{T}_u = \{1 \quad 0 \quad 0 \quad 0\} \tag{50a}$$

$$\mathbf{T}_w = \{0 \quad -1 \quad 0 \quad 0\} \tag{50b}$$

$$\mathbf{T}_\theta = \{0 \quad 0 \quad 1 \quad 0\} \tag{50c}$$

Subsequently, the expressions for each generalized displacement can be given as follows.

(1) Axial displacement of the beam axis

By integrating $u_{,x}(x)$, the expression of axial displacement can be expressed as

$$u(x) = u^a + \mathbf{N}_u(x)\mathbf{c} + U_u(x) \tag{51}$$

where $u^a$ refers to the axial displacement of the beam axis at the starting node, and $\mathbf{N}_u(x)$ and $U_u(x)$ can be obtained by

$$\mathbf{N}_u(x) = \mathbf{T}_u\mathbf{F}_\mathrm{n}\int_0^x \mathbf{N}_\sigma(\xi)\mathrm{d}\xi \tag{52}$$

$$U_u(x) = \mathbf{T}_u\mathbf{F}_\mathrm{n}\int_0^x \mathbf{F}_\sigma(\xi)\mathrm{d}\xi \tag{53}$$

(2) First-order differentiation of transverse displacement of the beam axis

By integrating $w_{,xx}(x)$, the expression of $w_{,x}(x)$ can be expressed as

$$w_{,x}(x) = w_{,x}^a + \mathbf{N}_{ww}(x)\mathbf{c} + U_{ww}(x) \tag{54}$$

where $w_{,x}^a$ is the first-order differentiation of transverse displacement at the starting node, and $\mathbf{N}_{ww}(x)$ and $U_{ww}(x)$ can be obtained by

$$\mathbf{N}_{ww}(x) = \mathbf{T}_w\mathbf{F}_\mathrm{n}\int_0^x \mathbf{N}_\sigma(\xi)\mathrm{d}\xi \tag{55}$$

$$U_{ww}(x) = \mathbf{T}_w\mathbf{F}_\mathrm{n}\int_0^x \mathbf{F}_\sigma(\xi)\mathrm{d}\xi \tag{56}$$

(3) Transverse displacement of the beam axis

By integrating $w_{,x}(x)$, the expression of $w(x)$ can be expressed as

$$w(x) = w^a + w_{,x}^a x + \mathbf{N}_w(x)\mathbf{c} + U_w(x) \tag{57}$$

where $w^a$ is the transverse displacement of the beam axis at the starting node, and $\mathbf{N}_w(x)$ and $U_w(x)$ can be obtained by

$$\mathbf{N}_w(x)=\mathbf{T}_w\mathbf{F}_{\mathrm{n}}\int_0^x\int_0^\xi \mathbf{N}_\sigma(\eta)\,\mathrm{d}\eta\mathrm{d}\xi \tag{58}$$

$$U_w(x)=\mathbf{T}_w\mathbf{F}_{\mathrm{n}}\int_0^x\int_0^\xi \mathbf{F}_\sigma(\eta)\,\mathrm{d}\eta\mathrm{d}\xi \tag{59}$$

(4) The magnitude of higher-order axial displacement

By integrating $\theta_{,x}(x)$, the magnitude of higher-order axial displacement can be expressed as

$$\theta(x)=\theta^a+\mathbf{N}_\theta(x)\mathbf{c}+U_\theta(x) \tag{60}$$

where $\theta^a$ represents the magnitude of higher-order axial displacement at the starting node, and $\mathbf{N}_\theta(x)$ and $U_\theta(x)$ can be obtained by

$$\mathbf{N}_\theta(x)=\mathbf{T}_\theta\mathbf{F}_{\mathrm{n}}\int_0^x \mathbf{N}_\sigma(\xi)\,\mathrm{d}\xi \tag{61}$$

$$U_\theta(x)=\mathbf{T}_\theta\mathbf{F}_{\mathrm{n}}\int_0^x \mathbf{F}_\sigma(\xi)\,\mathrm{d}\xi \tag{62}$$

(5) The measurement of the additional transverse displacement

By integrating $\phi_{,x}(x)$ as shown in Eq. (49) wherein $\theta(x)$ expressed in Eq. (60), the measurement of the additional transverse displacement can be expressed as

$$\phi(x)=\phi^a-\theta^a x+\mathbf{N}_\phi(x)\mathbf{c}+U_\phi(x) \tag{63}$$

where $\phi^a$ refers to the measurement of the additional transverse displacement at the starting node, and $\mathbf{N}_\phi(x)$ and $U_\phi(x)$ can be obtained by

$$\mathbf{N}_\phi(x)=D_{55}^{-1}\int_0^x \mathbf{N}_Q(\xi)\,\mathrm{d}\xi-\mathbf{T}_\theta\mathbf{F}_{\mathrm{n}}\int_0^x\int_0^\xi \mathbf{N}_\sigma(\eta)\,\mathrm{d}\eta\mathrm{d}\xi \tag{64}$$

$$U_\phi(x)=D_{55}^{-1}\int_0^x \mathbf{F}_Q(\xi)\,\mathrm{d}\xi-\mathbf{T}_\theta\mathbf{F}_{\mathrm{n}}\int_0^x\int_0^\xi \mathbf{F}_\sigma(\eta)\,\mathrm{d}\eta\mathrm{d}\xi \tag{65}$$

From the aforementioned expressions of the generalized displacements, it can be observed that the generalized displacement fields are functions of twelve parameters: the generalized displacements at the starting node ($u^a$, $w^a$, $w_{,x}^a$, $\theta^a$ and $\phi^a$) and the internal force parameters ($c_0 \sim c_6$). Notably, this parameter count precisely matches the total differential order of Eq. (18). For the sake of convenience, the generalized displacements can be integrated and represented as

$$\begin{Bmatrix} u(x) \\ w(x) \\ \theta(x) \\ \phi(x) \end{Bmatrix}=\begin{bmatrix} \mathbf{X}_u(x) \\ \mathbf{X}_w(x) \\ \mathbf{X}_\theta(x) \\ \mathbf{X}_\phi(x) \end{bmatrix}\mathbf{C}+\begin{Bmatrix} U_u(x) \\ U_w(x) \\ U_\theta(x) \\ U_\phi(x) \end{Bmatrix} \tag{66}$$

where

$$\mathbf{C}=\left\{u^a \quad w^a \quad w_{,x}^a \quad \theta^a \quad \phi^a \quad \mathbf{c}^{\mathrm{T}}\right\}^{\mathrm{T}} \tag{67}$$

$$\mathbf{X}_u(x)=\begin{bmatrix}1 & 0 & 0 & 0 & 0 & \mathbf{N}_u(x)\end{bmatrix} \tag{68a}$$

$$\mathbf{X}_w(x)=\begin{bmatrix}0 & 1 & x & 0 & 0 & \mathbf{N}_w(x)\end{bmatrix} \tag{68b}$$

$$\mathbf{X}_\theta(x)=\begin{bmatrix}0 & 0 & 0 & 1 & 0 & \mathbf{N}_\theta(x)\end{bmatrix} \tag{68c}$$

$$\mathbf{X}_{\phi}(x)=\begin{bmatrix}0 & 0 & 0 & -x & 1 & \mathbf{N}_{\phi}(x)\end{bmatrix} \tag{68d}$$

*3.3. Exact shape functions for generalized displacements*

For the development of a two-node exact beam element, each node should possess 6 DoFs to ensure that the total number of element displacement DoFs equals twelve. Consequently, $w_{,x}(x)$ and $\phi_{,x}(x)$ are added to the generalized displacement vector, which is expressed as

$$\mathbf{d}(x)=\mathbf{X}_{\mathrm{d}}(x)\mathbf{C}+\mathbf{U}_{\mathrm{d}}(x) \tag{69}$$

where

$$\mathbf{d}(x)=\{u(x) \quad w(x) \quad w_{,x}(x) \quad \theta(x) \quad \phi(x) \quad \phi_{,x}(x)\}^{\mathrm{T}} \tag{70}$$

$$\mathbf{X}_{\mathrm{d}}(x)=\begin{bmatrix}\mathbf{X}_{u}^{\mathrm{T}}(x) & \mathbf{X}_{w}^{\mathrm{T}}(x) & \mathbf{X}_{w,x}^{\mathrm{T}}(x) & \mathbf{X}_{\theta}^{\mathrm{T}}(x) & \mathbf{X}_{\phi}^{\mathrm{T}}(x) & \mathbf{X}_{\phi,x}^{\mathrm{T}}(x)\end{bmatrix}^{\mathrm{T}} \tag{71}$$

$$\mathbf{U}_{\mathrm{d}}(x)=\{U_{u}(x) \quad U_{w}(x) \quad U_{w,x}(x) \quad U_{\theta}(x) \quad U_{\phi}(x) \quad U_{\phi,x}(x)\}^{\mathrm{T}} \tag{72}$$

Meanwhile, the nodal displacement vector including the generalized displacements at two nodes of the element is defined by

$$\mathbf{\Phi}=\begin{bmatrix}\{u^{a} \quad w^{a} \quad w_{,x}^{a} \quad \theta^{a} \quad \phi^{a} \quad \phi_{,x}^{a}\}^{\mathrm{T}} \\ \{u^{b} \quad w^{b} \quad w_{,x}^{b} \quad \theta^{b} \quad \phi^{b} \quad \phi_{,x}^{b}\}^{\mathrm{T}}\end{bmatrix} \tag{73}$$

where the superscripts '$a$' and '$b$' correspond to the starting and ending nodes, respectively.

By substituting coordinates at the element ends ( $x=0,L$ ) into Eq. (69), the nodal generalized displacement vector can be expressed with the following expression according to the consistency condition that $\mathbf{\Phi}=\{\mathbf{d}^{\mathrm{T}}(0) \quad \mathbf{d}^{\mathrm{T}}(L)\}^{\mathrm{T}}$

$$\mathbf{\Phi}=\mathbf{X}\mathbf{C}+\mathbf{U} \tag{74}$$

where

$$\mathbf{X}=\begin{bmatrix}\mathbf{X}_{\mathrm{d}}(0) \\ \mathbf{X}_{\mathrm{d}}(L)\end{bmatrix},\ \mathbf{U}=\begin{bmatrix}\mathbf{U}_{\mathrm{d}}(0) \\ \mathbf{U}_{\mathrm{d}}(L)\end{bmatrix} \tag{75}$$

According to Eq. (74), the parameter vector $\mathbf{C}$ can be expressed by the nodal generalized displacement vector $\mathbf{\Phi}$ as

$$\mathbf{C}=\mathbf{X}^{-1}\mathbf{\Phi}-\mathbf{X}^{-1}\mathbf{U} \tag{76}$$

By substituting Eq. (76) into Eq. (66), the interpolation expression for the generalized displacements ($u(x)$, $w(x)$, $\theta(x)$ and $\phi(x)$) can be obtained as

$$u(x)=\mathbf{N}_{u}^{h}(x)\mathbf{\Phi}+U_{u}^{h}(x) \tag{77a}$$

$$w(x)=\mathbf{N}_{w}^{h}(x)\mathbf{\Phi}+U_{w}^{h}(x) \tag{77b}$$

$$\theta(x)=\mathbf{N}_{\theta}^{h}(x)\mathbf{\Phi}+U_{\theta}^{h}(x) \tag{77c}$$

$$\phi(x)=\mathbf{N}_{\phi}^{h}(x)\mathbf{\Phi}+U_{\phi}^{h}(x) \tag{77d}$$

where $\mathbf{N}_{u}^{h}(x)$, $\mathbf{N}_{w}^{h}(x)$, $\mathbf{N}_{\theta}^{h}(x)$ and $\mathbf{N}_{\phi}^{h}(x)$ are the exact shape function matrices, which are expressed as

$$\mathbf{N}_{u}^{h}(x)=\mathbf{X}_{u}(x)\mathbf{X}^{-1} \tag{78a}$$

$$\mathbf{N}_{w}^{h}(x)=\mathbf{X}_{w}(x)\mathbf{X}^{-1} \tag{78b}$$

$$\mathbf{N}_{\theta}^{h}(x)=\mathbf{X}_{\theta}(x)\mathbf{X}^{-1} \tag{78c}$$

$$\mathbf{N}_{\phi}^{h}(x)=\mathbf{X}_{\phi}(x)\mathbf{X}^{-1} \tag{78d}$$

and $U_{u}^{h}(x)$, $U_{w}^{h}(x)$, $U_{\theta}^{h}(x)$ and $U_{\phi}^{h}(x)$ represent the additional displacement distributions induced by the elemental load $q(x)$, which are expressed as

$$U_{u}^{h}(x)=U_{u}(x)-\mathbf{X}_{u}(x)\mathbf{X}^{-1}\mathbf{U} \tag{79a}$$

$$U_{w}^{h}(x)=U_{w}(x)-\mathbf{X}_{w}(x)\mathbf{X}^{-1}\mathbf{U} \tag{79b}$$

$$U_{\theta}^{h}(x)=U_{\theta}(x)-\mathbf{X}_{\theta}(x)\mathbf{X}^{-1}\mathbf{U} \tag{79c}$$

$$U_{\phi}^{h}(x)=U_{\phi}(x)-\mathbf{X}_{\phi}(x)\mathbf{X}^{-1}\mathbf{U} \tag{79d}$$

*3.4. Element stiffness equation and exact stiffness matrix*

Since the interpolation expressions for all generalized displacements have been obtained, the element stiffness equation can be derived according to the principle of virtual work. By introducing the constitutive equations of the cross-section as shown in Eqs. (14) and (15), Eq. (10) can be rewritten as

$$\int_{0}^{L}\left\{\delta\boldsymbol{\varepsilon}_{\mathrm{n}}^{\mathrm{T}}(x)\mathbf{D}_{\mathrm{n}}\boldsymbol{\varepsilon}_{\mathrm{n}}(x)+\delta\left[\theta(x)+\phi_{,x}(x)\right]D_{55}\left[\theta(x)+\phi_{,x}(x)\right]\right\}\mathrm{d}x-\int_{0}^{L}\delta w(x)q(x)\mathrm{d}x=\delta\mathbf{\Phi}^{\mathrm{T}}\mathbf{P}_{e} \tag{80}$$

where

$$\mathbf{P}_{e}=\left\{P_{ax}\quad P_{ay}\quad M_{a}\quad 0\quad 0\quad 0\quad P_{bx}\quad P_{by}\quad M_{b}\quad 0\quad 0\quad 0\right\}^{\mathrm{T}} \tag{81}$$

Based on the interpolation expressions of generalized displacements, the generalized strains and their variation expressions can be obtained by

$$\boldsymbol{\varepsilon}_{\mathrm{n}}(x)=\overline{\mathbf{N}}_{\mathrm{n}}(x)\mathbf{\Phi}+\overline{\mathbf{U}}_{\mathrm{n}}(x) \tag{82}$$

$$\delta\boldsymbol{\varepsilon}_{\mathrm{n}}(x)=\overline{\mathbf{N}}_{\mathrm{n}}(x)\delta\mathbf{\Phi} \tag{83}$$

$$\theta(x)+\phi_{,x}(x)=\overline{\mathbf{N}}_{\mathrm{s}}(x)\mathbf{\Phi}+\overline{U}_{\mathrm{s}}(x) \tag{84}$$

$$\delta\left[\theta(x)+\phi_{,x}(x)\right]=\overline{\mathbf{N}}_{\mathrm{s}}(x)\delta\mathbf{\Phi} \tag{85}$$

where

$$\overline{\mathbf{N}}_{\mathrm{n}}(x)=\begin{bmatrix}\mathbf{N}_{u,x}^{h}(x)\\-\mathbf{N}_{w,xx}^{h}(x)\\\mathbf{N}_{\theta,x}^{h}(x)\\\mathbf{N}_{\phi}^{h}(x)\end{bmatrix} \tag{86}$$

$$\overline{\mathbf{U}}_{\mathrm{n}}(x)=\begin{bmatrix}U_{u.x}^{h}(x)\\-U_{w,xx}^{h}(x)\\U_{\theta,x}^{h}(x)\\U_{\phi}^{h}(x)\end{bmatrix} \tag{87}$$

$$\overline{\mathbf{N}}_{\mathrm{s}}(x)=\mathbf{N}_{\theta}^{h}(x)+\mathbf{N}_{\phi,x}^{h}(x) \tag{88}$$

$$\overline{U}_{\mathrm{s}}(x)=U_{\theta}^{h}(x)+U_{\phi,x}^{h}(x) \tag{89}$$

By substituting Eqs. (82)-(85) into Eq. (80) and using the variation of $\delta w(x)$ as $\delta w(x)=\mathbf{N}_{w}^{h}(x)\delta\mathbf{\Phi}$, the element stiffness equations can be written as

$$\mathbf{K}_e \mathbf{\Phi} = \mathbf{F}_e + \mathbf{P}_e \tag{90}$$

where $\mathbf{K}_e$ and $\mathbf{F}_e$ represent the exact element stiffness matrix and equivalent nodal force vector, respectively, they are expressed as

$$\mathbf{K}_e = \int_0^L \left\{ \bar{\mathbf{N}}_{\mathrm{n}}^{\mathrm{T}}(x) \mathbf{D}_{\mathrm{n}} \bar{\mathbf{N}}_{\mathrm{n}}(x) \mathbf{\Phi} + \bar{\mathbf{N}}_{\mathrm{s}}^{\mathrm{T}}(x) D_{55} \bar{\mathbf{N}}_{\mathrm{s}}(x) \mathbf{\Phi} \right\} \mathrm{d}x \tag{91}$$

$$\mathbf{F}_e = \int_0^L \left[ \mathbf{N}_w^h(x) \right]^{\mathrm{T}} q(x) \mathrm{d}x - \int_0^L \left\{ \bar{\mathbf{N}}_{\mathrm{n}}^{\mathrm{T}}(x) \mathbf{D}_{\mathrm{n}} \bar{\mathbf{U}}_{\mathrm{n}}(x) + \bar{\mathbf{N}}_{\mathrm{s}}^{\mathrm{T}}(x) D_{55} \bar{U}_{\mathrm{s}}(x) \right\} \mathrm{d}x \tag{92}$$

Then, the stiffness equations of the structure can be obtained by using the conventional assembly methods.

## 4. Improvement of stress expressions

### *4.1. Equilibrium-based stress expressions*

The aforementioned derivation of the element stiffness matrix and element stiffness equation based on exact shape functions effectively addresses the issue of finite element discretization errors. However, the stress solutions cannot strictly satisfy the differential equilibrium equations, which limits the accuracy of quasi-3D beam models. Regarding the problem of stress distribution deviations, the transverse shear stress and transverse normal stress should be derived by further incorporating the differential equilibrium equations that govern the relationships between stresses, which are expressed as

$$\frac{\partial \sigma_x(x,y)}{\partial x} + \frac{\partial \tau_{xy}(x,y)}{\partial y} = 0 \tag{93}$$

$$\frac{\partial \sigma_y(x,y)}{\partial y} + \frac{\partial \tau_{xy}(x,y)}{\partial x} = 0 \tag{94}$$

where the body forces along $x$-axis and $y$-axis are neglected for the sake of formulation simplicity.

Firstly, according to Eqs. (5), (2), (3) and (22), the expression of axial normal stress $\sigma_x(x,y)$ can be rewritten as

$$\sigma_x(x,y) = \mathbf{T}_x(y) \boldsymbol{\varepsilon}_{\mathrm{n}}(x) = \mathbf{T}_x(y) \mathbf{F}_{\mathrm{n}} \boldsymbol{\sigma}_{\mathrm{n}}(x) \tag{95}$$

where

$$\mathbf{T}_x(y) = \left\{ \bar{Q}_{11}(y) \quad y\bar{Q}_{11}(y) \quad f(y)\bar{Q}_{11}(y) \quad g_{,y}(y)\bar{Q}_{12}(y) \right\} \tag{96}$$

Using Eqs. (93) and (95) and neglecting the impact of axial force and transverse stretching force on transverse shear stress, the equilibrium-based transverse shear stress can be expressed as [23]

$$\tau_{xy}(x,y) = \mathbf{T}_{xy}(y) \boldsymbol{\tau}(x) \tag{97}$$

where

$$\boldsymbol{\tau}(x) = \left\{ M_{w,x}(y) \quad M_{\theta,x}(x) \right\}^{\mathrm{T}} \tag{98}$$

$$\mathbf{T}_{xy}(y) = \left\{ S_w(y) \quad S_\theta(y) \right\} \tag{99}$$

with

$$S_w(y) = -\int_{-\frac{h}{2}}^{y} \mathbf{T}_x(\xi) \left[ f_{12} \quad f_{22} \quad f_{23} \quad f_{24} \right]^{\mathrm{T}} \mathrm{d}\xi \tag{100a}$$

$$S_\theta(y) = -\int_{-\frac{h}{2}}^{y} \mathbf{T}_x(\xi) \left[ f_{13} \quad f_{23} \quad f_{33} \quad f_{34} \right]^{\mathrm{T}} \mathrm{d}\xi \tag{100b}$$

By substituting Eq. (97) into Eq. (94), the equilibrium-based transverse normal stress can be expressed as

$$\sigma_y(x,y) = c_\sigma - \int_{-\frac{h}{2}}^{y} \frac{\partial \tau_{xy}(x,\xi)}{\partial x} \mathrm{d}\xi = \mathbf{T}_y(y)\boldsymbol{\tau}_{,x}(x) \tag{101}$$

where the integration constant $c_\sigma$ is zero because the transverse normal stress at $y = -h/2$ is zero, and $\mathbf{T}_{\bar{y}}(y)$ can be expressed as

$$\mathbf{T}_y(y) = \left\{ -\int_{-\frac{h}{2}}^{y} S_w(\xi)\mathrm{d}\xi \quad -\int_{-\frac{h}{2}}^{y} S_\theta(\xi)\mathrm{d}\xi \right\} \tag{102}$$

*4.2. Modified cross-sectional stiffness*

The equilibrium-based stress formulations, differing substantially from conventional expressions, necessitate modified cross-sectional constitutive equations. Following the framework of Ref. [20], we apply the energy equivalence principle to develop the equilibrium-based constitutive models related to normal stress and shear stress, respectively.

The energy work by equilibrium-based normal stresses can be expressed as

$$\begin{aligned}\Pi_\sigma = &\int_0^L \int_A \left[ \left\{ \sigma_x(x,y) \quad \sigma_y(x,y) \right\} \begin{Bmatrix} \varepsilon_x(x,y) \\ \varepsilon_y(x,y) \end{Bmatrix} \right] \mathrm{d}A\mathrm{d}x \\ &-\frac{1}{2}\int_0^L \int_A \left[ \left\{ \sigma_x(x,y) \quad \sigma_y(x,y) \right\} \begin{bmatrix} s_{11}(y) & s_{12}(y) \\ s_{21}(y) & s_{22}(y) \end{bmatrix} \begin{Bmatrix} \sigma_x(x,y) \\ \sigma_y(x,y) \end{Bmatrix} \right] \mathrm{d}A\mathrm{d}x \end{aligned} \tag{103}$$

where

$$\begin{bmatrix} s_{11}(y) & s_{12}(y) \\ s_{21}(y) & s_{22}(y) \end{bmatrix} = \begin{bmatrix} \bar{Q}_{11}(y) & \bar{Q}_{12}(y) \\ \bar{Q}_{12}(y) & \bar{Q}_{11}(y) \end{bmatrix}^{-1} \tag{104}$$

By substituting Eqs. (95), (101) and (2), Eq. (103) can be further rewritten as

$$\Pi_\sigma = \int_0^L \left\{ \begin{bmatrix} \boldsymbol{\varepsilon}_{\mathrm{n}}(x) \\ \boldsymbol{\tau}_{,x}(x) \end{bmatrix}^{\mathrm{T}} \begin{bmatrix} \mathbf{H}_{xx}\boldsymbol{\varepsilon}_{\mathrm{n}}(x) \\ \mathbf{H}_{yy}\boldsymbol{\varepsilon}_{\mathrm{n}}(x) \end{bmatrix} - \frac{1}{2} \begin{bmatrix} \boldsymbol{\varepsilon}_{\mathrm{n}}(x) \\ \boldsymbol{\tau}_{,x}(x) \end{bmatrix}^{\mathrm{T}} \begin{bmatrix} \mathbf{H}_{xsx} & \mathbf{H}_{xsy} \\ \mathbf{H}_{xsy}^{\mathrm{T}} & \mathbf{H}_{ysy} \end{bmatrix} \begin{bmatrix} \boldsymbol{\varepsilon}_{\mathrm{n}}(x) \\ \boldsymbol{\tau}_{,x}(x) \end{bmatrix} \right\} \mathrm{d}x \tag{105}$$

where

$$\mathbf{H}_{xx} = \int_A \mathbf{T}_x^{\mathrm{T}}(y)\mathbf{B}_x(y)\mathrm{d}A \tag{106a}$$

$$\mathbf{H}_{yy} = \int_A \mathbf{T}_y^{\mathrm{T}}(y)\mathbf{B}_y(y)\mathrm{d}A \tag{106b}$$

$$\mathbf{H}_{xsx} = \int_A \mathbf{T}_x^{\mathrm{T}}(y)s_{11}(y)\mathbf{T}_x(y)\mathrm{d}A \tag{106c}$$

$$\mathbf{H}_{ysy} = \int_A \mathbf{T}_y^{\mathrm{T}}(y)s_{22}(y)\mathbf{T}_y(y)\mathrm{d}A \tag{106d}$$

$$\mathbf{H}_{xsy} = \int_A \mathbf{T}_x^{\mathrm{T}}(y)s_{12}(y)\mathbf{T}_y(y)\mathrm{d}A \tag{106e}$$

By postulating the existence of an internal force vector $\boldsymbol{\sigma}_{\mathrm{n}}(x) = \mathbf{D}_{\mathrm{n}}^{\mathrm{m}}\boldsymbol{\varepsilon}_{\mathrm{n}}(x)$ with $\mathbf{D}_{\mathrm{n}}^{\mathrm{m}}$ representing the modified cross-sectional stiffness matrix related to the normal stresses, the energy performed by these internal forces can be expressed as $U_\sigma = \frac{1}{2}\int_0^L \boldsymbol{\sigma}_{\mathrm{n}}^{\mathrm{T}}(\boldsymbol{\varepsilon}_{\mathrm{n}}(x))\boldsymbol{\varepsilon}_{\mathrm{n}}(x)\mathrm{d}x$ and the energy associated with equilibrium-based stresses should be equivalent to $U_\sigma$, namely, $\Pi_\sigma = U_\sigma$. As a sufficient condition to guarantee the validity of $\Pi_\sigma = U_\sigma$, the subsequent equation formulated based on the differential segment of the beam should hold

$$\begin{bmatrix} \boldsymbol{\varepsilon}_{\mathrm{n}}(x) \\ \boldsymbol{\tau}_{,x}(x) \end{bmatrix}^{\mathrm{T}} \begin{bmatrix} \mathbf{H}_{xx}\boldsymbol{\varepsilon}_{\mathrm{n}}(x) \\ \mathbf{H}_{yy}\boldsymbol{\varepsilon}_{\mathrm{n}}(x) \end{bmatrix} - \frac{1}{2} \begin{bmatrix} \boldsymbol{\varepsilon}_{\mathrm{n}}(x) \\ \boldsymbol{\tau}_{,x}(x) \end{bmatrix}^{\mathrm{T}} \begin{bmatrix} \mathbf{H}_{xsx} & \mathbf{H}_{xsy} \\ \mathbf{H}_{xsy}^{\mathrm{T}} & \mathbf{H}_{ysy} \end{bmatrix} \begin{bmatrix} \boldsymbol{\varepsilon}_{\mathrm{n}}(x) \\ \boldsymbol{\tau}_{,x}(x) \end{bmatrix} = \frac{1}{2}\boldsymbol{\sigma}_{\mathrm{n}}^{\mathrm{T}}(\boldsymbol{\varepsilon}_{\mathrm{n}}(x))\boldsymbol{\varepsilon}_{\mathrm{n}}(x) \tag{107}$$

Based on the variation principle and considering that $\delta\boldsymbol{\varepsilon}_{\mathrm{n}}(x)$ and $\delta\boldsymbol{\tau}_{,x}(x)$ are arbitrary variations, the following two sets of equations can be derived

$$\mathbf{H}_{\varepsilon x}\boldsymbol{\varepsilon}_{\mathrm{n}}(x)+\mathbf{H}_{\varepsilon y}^{\mathrm{T}}\boldsymbol{\tau}_{,x}(x)=\boldsymbol{\sigma}_{\mathrm{n}}(x) \tag{108}$$

$$\mathbf{H}_{\varepsilon y}\boldsymbol{\varepsilon}_{\mathrm{n}}(x)-\mathbf{H}_{ysy}\boldsymbol{\tau}_{,x}(x)=\mathbf{0} \tag{109}$$

where

$$\mathbf{H}_{\varepsilon x}=\mathbf{H}_{xx}+\mathbf{H}_{xx}^{\mathrm{T}}-\mathbf{H}_{xsx} \tag{110a}$$

$$\mathbf{H}_{\varepsilon y}=\mathbf{H}_{yy}-\frac{1}{2}\mathbf{H}_{xsy}-\frac{1}{2}\mathbf{H}_{xsy}^{\mathrm{T}} \tag{110b}$$

Then, the following equation reflecting the cross-sectional constitutive equations related to normal stresses can be obtained from Eqs. (108) and (109) as

$$\mathbf{D}_{\mathrm{n}}^{\mathrm{m}}\boldsymbol{\varepsilon}_{\mathrm{n}}(x)=\boldsymbol{\sigma}_{\mathrm{n}}(x) \tag{111}$$

where $\mathbf{D}_{\mathrm{n}}^{\mathrm{m}}$ can be expressed as

$$\mathbf{D}_{\mathrm{n}}^{\mathrm{m}}=\mathbf{H}_{\varepsilon x}+\mathbf{H}_{\varepsilon y}^{\mathrm{T}}\mathbf{H}_{ysy}^{-1}\mathbf{H}_{\varepsilon y} \tag{112}$$

For the shear stiffness, assuming the existence of the shear force  that satisfies the cross-sectional constitutive equation, the energy performed by this shear force can be expressed as $U_{\tau}=\frac{1}{2}\int_{0}^{L}Q(x)\left[\theta(x)+\phi_{,x}(x)\right]\mathrm{d}x$. Similar to the method above, the following constitutive equation for cross-sectional shear deformation can be derived based on the consistency between the above energy expression and the energy work by equilibrium-based shear stress given in Ref. [20]

$$D_{55}^{\mathrm{m}}\left[\theta(x)+\phi_{,x}(x)\right]=Q(x) \tag{113}$$

where $D_{55}^{\mathrm{m}}$ represents the modified cross-sectional shear stiffness, which is expressed as

$$D_{55}^{\mathrm{m}}=\mathbf{f}_{s}^{\mathrm{T}}\mathbf{f}_{ss}^{-1}\mathbf{f}_{s} \tag{114}$$

with

$$\mathbf{f}_{s}=\int_{A}\mathbf{T}_{xy}^{\mathrm{T}}(y)g(y)\mathrm{d}A \tag{115}$$

$$\mathbf{f}_{ss}=\int_{A}\left[\mathbf{T}_{xy}^{\mathrm{T}}(y)\mathbf{T}_{xy}(y)/G(y)\right]\mathrm{d}A \tag{116}$$

Thus, the equilibrium-based modified cross-sectional stiffness terms $\mathbf{D}_{\mathrm{n}}^{\mathrm{m}}$ and $D_{55}^{\mathrm{m}}$ have been respectively derived. Their substitution for $\mathbf{D}_{\mathrm{n}}$ and $D_{55}$ in **Sec. 3** directly yields the formulation of an equilibrium-based quasi-3D beam exact finite element.

## 5. Numerical examples

### *5.1. FG material models*

Two types of FG sandwich material models, denoted as Type A and Type B, respectively, are considered in the section of numerical examples. Type A represents the sandwich model with FG outer surfaces and a homogeneous core, and Type B refers to the sandwich model featuring FG core and homogeneous outer surfaces. For the sake of simplicity, it is assumed that Poisson’s ratio is constant. The Young’s modulus that varies along with the thickness is expressed as

$$E(y)=E_m+(E_c-E_m)V_c(y) \tag{117}$$

where $E(y)$ is Young's modulus at the location of $y$, $E_c$ and $E_m$ refers to the Young's moduli of ceramic material and metal material, respectively, and $V_c(y)$ represents the volume fraction. For Type A and Type B materials, $V_c(y)$ is defined as:

(1) Type A: sandwich model with FG surfaces and homogeneous cores

$$V_c(y)=\begin{cases}\left[(y-h_0)/(h_1-h_0)\right]^p & for\ \ y\in[h_0,h_1] \\ 1 & for\ \ y\in[h_1,h_2] \\ \left[(y-h_3)/(h_2-h_3)\right]^p & for\ \ y\in[h_2,h_3]\end{cases} \tag{118}$$

(2) Type B: sandwich model with FG core and homogeneous surfaces

$$V_c(y)=\begin{cases}0 & for\ \ y\in[h_0,h_1] \\ \left[(y-h_1)/(h_2-h_1)\right]^p & for\ \ y\in[h_1,h_2] \\ 1 & for\ \ y\in[h_2,h_3]\end{cases} \tag{119}$$

where $p$ is the power-law index and $y$ is the coordinate in thickness direction, $h_0,h_1,h_2$ and $h_3$ are characteristic positions related to material distribution, including the junction position of adjacent material layers and the boundary position of the beam's thickness, as shown in **Fig. 3**.

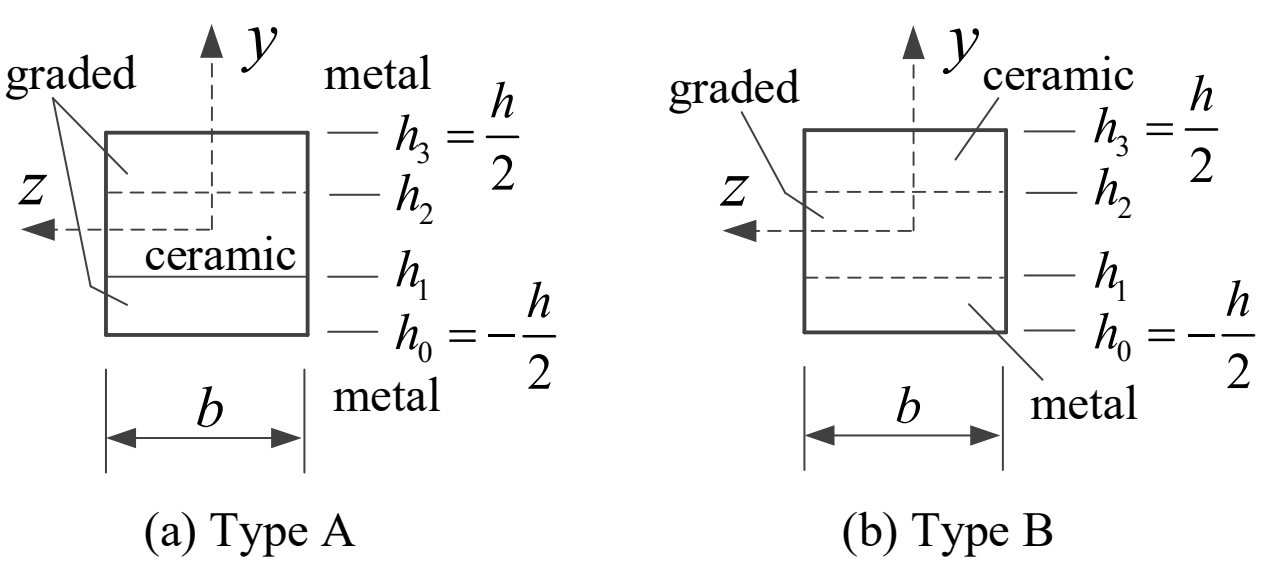

**Fig. 3.** The characteristic positions related to material distribution.

In the following numerical investigation, the FG material properties are set to be [30]: Aluminum ( $\mathrm{Al}: E_m = 70000\,\mathrm{N/mm^2}$ ) and Alumina ($Al_2O_3$: $\mathrm{Al}: E_c = 380000\,\mathrm{N/mm^2}$ ). The Poisson's ratio of the material is set to $v = 0.3$ . The characteristic positions are set to $h_0 = -100\mathrm{mm}, h_1 = -40\mathrm{mm}, h_2 = 40\mathrm{mm}$ and $h_3 = 100\mathrm{mm}$ .

*5.2. Verification of computational performance*

In this section, the FG sandwich beams under a uniformly distributed load applied at the upper surface are investigated. The load is set as $q(x) = -50\mathrm{N/mm}$ , the size of the cross-section is set to $b\times h = 50\mathrm{mm}\times 200\mathrm{mm}$ and the length of the beam is $L = 2000\mathrm{mm}$ . Specifically, two boundary conditions including Simply-Supported (SS) and Clamped-Clamped (CC) boundaries are considered, as shown in **Fig. 4**.

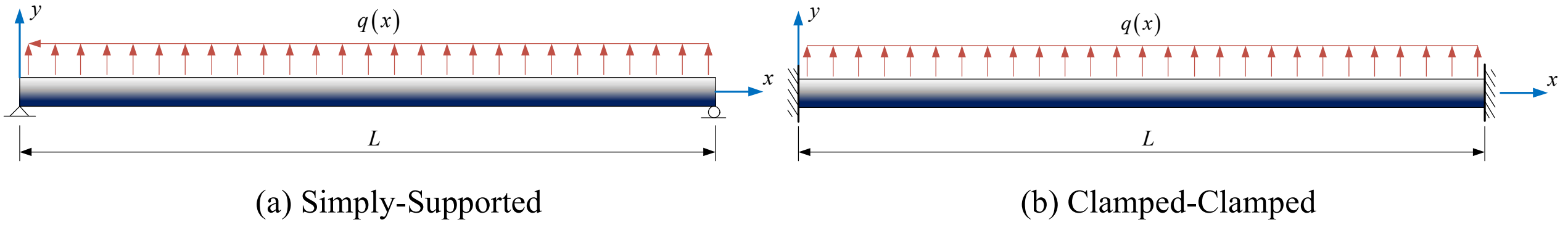

(a) Simply-Supported (b) Clamped-Clamped

**Fig. 4.** Two boundary conditions for the FG sandwich beam.

The following four beam element models are used to predict the displacement and stress solutions. If not otherwise specified, the third-order polynomial function (see **Table 1**) is employed to describe the higher-order distribution of out-of-plane axial displacement. The four element models are denoted as follows:

(1) DTS: the conventional Displacement-based beam finite element based on Third-order Shear deformation theory. The element has two nodes. For discretization, linear interpolation is used for the axial displacement of the beam axis and the magnitude of higher-order axial displacement, and cubic Hermite interpolation is performed for the transverse displacement of the beam axis. Then, the number of DoFs for each node is 4. As introduced in Ref. [23], the element stiffness matrix of DTS is formulated based on the principle of virtual work.

(2) DQ3D: the Displacement-based beam finite element based on Quasi-3D theory, which is a two-node beam element, where cubic Hermite interpolation is performed for all generalized displacements, including the axial and transverse displacements of the beam axis, the magnitude of higher-order axial displacement and the measurement of the additional transverse displacement from transverse stretching deformation. Then, the number of DoFs for each node is 8. In this element, the strains and stresses are computed by using Eq. (2) and Eq. (5), respectively. The element stiffness matrix is derived based on the principle of virtual work with generalized displacements as the unknown fields. Readers are referred to Ref. [30] for details of DQ3D, which is denoted as 'DTS-3D' in Ref. [30].

(3) MQ3D: the Mixed beam finite element based on quasi-3D theory, which is the two-node beam element proposed by Li et al. [30]. In this element, linear interpolation is used for discretization of the generalized displacement components $u(x)$, $\theta(x)$ and $\phi(x)$, and cubic Hermite interpolation is adopted for discretization of the generalized displacement component $w(x)$. Then, the number of DoFs for each node is 5. Especially, the differential equilibrium equation governing the relationship between $\sigma_x(x,y)$ and $\tau_{xy}(x,y)$ (Eq. (93)) is employed to derive the expression of transverse shear stress, while the normal stresses are obtained via conventional method. Readers are referred to Ref. [30] for details of MQ3D, which is denoted as 'MTS-3D' in Ref. [30].

(4) FEQ3D: the Force-based Exact beam finite element based on Quasi-3D theory, which is implemented according to the formulation presented in this work. The element has two nodes and the number of DoFs for each node is 6. The two differential equilibrium equations (Eqs. (93) and (94)) are introduced and the axial normal stress, transverse shear stress and transverse normal stress are obtained by using Eqs. (95), (97) and (101), respectively.

In addition, in contrast to the beam element previously listed, the 4-node planar Quadrilateral element (denoted as Q4) is capable of comprehensively accounting for diverse in-plane deformations. Note that the plane stress assumption is employed for the Q4 element in this study. Typically, precise displacement and stress outcomes can be achieved by the Q4 element when the mesh is refined adequately. Therefore, the solutions produced by Q4 are taken as the references for the evaluation of the beam elements. Specifically, the solutions of Q4 are obtained with a mesh of $m_x \times m_y = 400 \times 100$, where $m_x$ and $m_y$ denote the number of elements along the *x*-axis and the *y*-axis, respectively.

The number of Q4 elements through the thickness of the bottom layer, core layer and top layer is 30, 40 and 30, respectively, and this setting enables an accurate representation of the beam's deformation and stress distribution.

It is important to note that a well-established fact for FG beams is the deviation of the neutral surface from the geometric mid-surface [62, 63], which can lead to discrepancies between modeling approaches. The present model is formulated based on the geometric mid-surface to align with many classical models, allowing for a straightforward verification of the proposed element's performance. Readers should be aware that potential discrepancies between the present results and those from neutral-surface-based models are expected in cases where the neutral surface shift is pronounced.

Firstly, a convergence investigation is conducted for the four beam elements under the setting of $p = 5.0$, and the results are presented in **Table 2**-**Table 5** for different settings of boundary conditions and material models. Specifically, details for the solution of the SS beam with the Type B material model using FEQ3D are provided in **Appendix**. The results indicate that the proposed beam element (FEQ3D) works out convergent displacement solutions by using only one element, which verifies the accuracy of the closed-form solutions for the internal force fields and the rationality of the theoretical framework for constructing an exact finite element model based on the analytical internal force fields. For other elements (DTS, DQ3D and MQ3D), whether displacement-based or mixed formulations, the discretization error is inevitable due to their polynomial approximation nature, necessitating mesh refinement to achieve convergent results. In **Table 2**-**Table 5**, the total number of DoFs for the four elements to achieve converged solutions is also listed. The data demonstrate that DTS, DQ3D and MQ3D require a large number of DoFs to maintain computational stability, particularly for clamped-clamped boundary conditions. In contrast, the proposed element achieves convergent results with only 12 DoFs. In the following analyses for displacements and stresses, numerical modeling is carried out with the meshes required for convergent solutions in the convergence investigation.

**Table 2**-**Table 5** further present a comparative analysis of the computational time required to achieve convergent displacement solutions, where $T_{\mathrm{DTS}}$ and $T_{\mathrm{beam}}$ denote the computational time consumed by DTS and other beam element models, respectively. The results indicate that for simply-supported beams, which exhibit relatively simple deformation and internal force distributions, DQ3D, MQ3D and FEQ3D demonstrate high computational efficiency. However, in the case of clamped-clamped beams with more complex mechanical behavior, the proposed FEQ3D significantly outperforms the other three models (DTS, DQ3D and MQ3D), highlighting its superior computational efficiency.

For the two material models and two boundary conditions, the convergent vertical displacement solutions at mid-span obtained by the four beam elements and the Q4 model are listed in **Table 6** and **Table 7**, under different settings of power-law index. For better comparison, the Relative Errors (RE) under various power-law indices and the Average Relative Errors (ARE) for the four beam elements are provided in the tables.

As shown by the results in **Table 6** and **Table 7**, FEQ3D has significantly higher solution accuracy than the other beam elements. The superior performance stems from the rational incorporation of both equilibrium-based transverse shear and normal stresses. The conventional quasi-3D beam element (DQ3D) exhibits significant displacement solution errors due to inaccuracies in its stress expressions, particularly for asymmetric material distributions (Type B). While MQ3D improves upon DQ3D by introducing equilibrium-based transverse shear stress through the mixed formulation method, its displacement accuracy remains limited by the persistent inaccuracies in transverse normal stress representation. This evidence conclusively demonstrates that equilibrium-based transverse normal stress plays a pivotal role in enhancing displacement solution accuracy for FG sandwich beams. The DTS element achieves even higher displacement solution precision than both DQ3D and MQ3D by omitting transverse stretching deformation and neglecting transverse normal stress effects. This strategic simplification effectively circumvents the errors that would otherwise arise from the inherent inaccuracies in the transverse normal stress expressions. Nevertheless, DTS's displacement solution accuracy remains inferior to that of FEQ3D, and its inability to predict transverse normal stress significantly limits its practical utility despite the advantages as mentioned above.

**Table 2** Convergence of the mid-span vertical displacement [$10^{-2}$ mm] (SS, Type A, $p = 5.0$ ).

| Number of elements | DTS [23] | DQ3D [30] | MQ3D [30] | FEQ3D |
|---|---|---|---|---|
| 1 | | | | 280.48 |
| 2 | 170.22 | 280.52 | 280.56 | 280.48 |
| 4 | 253.70 | 280.45 | 280.58 | 280.48 |
| 8 | 274.52 | 280.42 | 280.61 | |
| 16 | 279.63 | 280.41 | 280.61 | |
| 32 | 280.80 | 280.41 | | |
| 64 | 281.05 | | | |
| 128 | 281.10 | | | |
| 256 | 281.12 | | | |
| 512 | 281.12 | | | |
| Convergent solution | 281.12 | 280.41 | 280.61 | 280.48 |
| Number of DoFs | 1285 | 136 | 45 | 12 |
| $T_{\text{beam}}/T_{\text{DTS}}$ | 1.000 | 0.214 | 0.188 | 0.291 |

**Table 3** Convergence of the mid-span vertical displacement [$10^{-2}$ mm] (SS, Type B, $p = 5.0$ ).

| Number of elements | DTS [23] | DQ3D [30] | MQ3D [30] | FEQ3D |
|---|---|---|---|---|
| 1 | | | | 232.68 |
| 2 | 142.09 | 223.07 | 225.39 | 232.68 |
| 4 | 209.59 | 222.87 | 224.55 | 232.68 |
| 8 | 225.94 | 222.82 | 224.24 | |
| 16 | 229.62 | 222.82 | 224.16 | |
| 32 | 230.43 | | 224.14 | |
| 64 | 230.63 | | 224.14 | |
| 128 | 230.68 | | | |
| 256 | 230.69 | | | |
| 512 | 230.69 | | | |
| Convergent solution | 230.69 | 222.82 | 224.14 | 232.68 |
| Number of DoFs | 1285 | 72 | 165 | 12 |
| $T_{\text{beam}}/T_{\text{DTS}}$ | 1.000 | 0.157 | 0.328 | 0.291 |

The stress solutions are further investigated. For the setting of $p = 5.0$, the distributions of stresses through the thickness at $x = 1500\text{mm}$, obtained by different element models, are presented in **Fig. 5**-**Fig. 10**, for different stress components under the two boundary conditions.

**Table 4** Convergence of the mid-span vertical displacement $[10^{-2}$ mm] (CC, Type A, $p = 5.0$).

| Number of elements | DTS [23] | DQ3D [30] | MQ3D [30] | FEQ3D |
|---|---|---|---|---|
| 1 | | | | 58.003 |
| 2 | 3.2419 | 56.001 | 58.676 | 58.003 |
| 4 | 44.973 | 56.659 | 58.828 | 58.003 |
| 8 | 55.359 | 57.057 | 58.775 | |
| 16 | 57.842 | 57.449 | 58.550 | |
| 32 | 58.337 | 57.746 | 58.488 | |
| 64 | 58.418 | 57.794 | 58.472 | |
| 128 | 58.437 | 57.796 | 58.468 | |
| 256 | 58.442 | 57.796 | 58.467 | |
| 512 | 58.444 | | 58.466 | |
| 1024 | 58.444 | | 58.466 | |
| Convergent solution | 58.444 | 57.796 | 58.466 | 58.003 |
| Number of DoFs | 2565 | 1032 | 2565 | 12 |
| $T_{\text{beam}}/T_{\text{DTS}}$ | 1.000 | 0.288 | 1.737 | 0.086 |

**Table 5** Convergence of the mid-span vertical displacement $[10^{-2}$ mm] (CC, Type B, $p = 5.0$).

| Number of elements | DTS [23] | DQ3D [30] | MQ3D [30] | FEQ3D |
|---|---|---|---|---|
| 1 | | | | 50.880 |
| 2 | 6.5794 | 46.760 | 49.806 | 50.880 |
| 4 | 40.237 | 47.300 | 49.854 | 50.880 |
| 8 | 48.085 | 47.639 | 49.769 | |
| 16 | 49.549 | 47.880 | 49.595 | |
| 32 | 49.815 | 48.039 | 49.521 | |
| 64 | 49.894 | 48.069 | 49.508 | |
| 128 | 49.920 | 48.071 | 49.504 | |
| 256 | 49.927 | 48.071 | 49.503 | |
| 512 | 49.929 | | 49.503 | |
| 1024 | 49.929 | | | |
| Convergent solution | 49.929 | 48.071 | 49.503 | 50.880 |
| Number of DoFs | 2565 | 1032 | 1285 | 12 |
| $T_{\text{beam}}/T_{\text{DTS}}$ | 1.000 | 0.288 | 0.525 | 0.086 |

As illustrated in **Fig. 5** and **Fig. 8**, the axial normal stress distribution results obtained from DQ3D, MQ3D and FEQ3D exhibit a high degree of consistency with those derived from Q4, for both material models. This indicates that all three beam models based on quasi-3D theory can accurately predict axial stress outcomes. As evident from **Fig. 7** and **Fig. 10**, the transverse shear stress distributions obtained from DQ3D exhibit significant discrepancies compared to the results from Q4, particularly at the material interfaces ($y = \pm 40\text{mm}$). The results indicate that the conventional approach of directly deriving transverse shear stress solely from kinematic assumptions and constitutive equations is inadequate, even when the higher-order function has been incorporated into the element formulation to account for non-uniform shear stress effects. For MQ3D and FEQ3D, since they employ equilibrium-based transverse shear stress formulations, their predicted transverse shear stress distributions exhibit good agreement with the Q4 results. As illustrated, the accurate transverse shear stress should vary as a continuous and smooth curve through the thickness direction.

**Table 6** Comparison of the mid-span vertical displacement solutions [$10^{-2}$ mm] (SS).

| Type | *p* | DTS [23] | | DQ3D [30] | | MQ3D [30] | | FEQ3D | | Q4 |
|---|---|---|---|---|---|---|---|---|---|---|
| | | Disp | RE (%) | Disp | RE (%) | Disp | RE (%) | Disp | RE (%) | |
| A | 0.0 | 84.289 | 0.011 | 84.110 | 0.202 | 84.110 | 0.202 | 84.109 | 0.203 | 84.280 |
| | 0.5 | 126.59 | 0.182 | 126.29 | 0.052 | 126.30 | 0.047 | 126.30 | 0.047 | 126.36 |
| | 1.0 | 162.00 | 0.229 | 161.61 | 0.012 | 161.61 | 0.012 | 161.61 | 0.012 | 161.63 |
| | 2.0 | 212.91 | 0.302 | 212.37 | 0.048 | 212.40 | 0.061 | 212.39 | 0.057 | 212.27 |
| | 5.0 | 281.12 | 0.300 | 280.41 | 0.047 | 280.61 | 0.118 | 280.48 | 0.071 | 280.28 |
| | 10.0 | 314.74 | 0.220 | 313.96 | 0.029 | 314.41 | 0.115 | 314.07 | 0.006 | 314.05 |
| B | 0.0 | 166.35 | 0.078 | 163.27 | 1.928 | 163.51 | 1.784 | 166.80 | 0.192 | 166.48 |
| | 0.5 | 192.63 | 0.109 | 187.63 | 2.705 | 187.83 | 2.598 | 193.57 | 0.379 | 192.84 |
| | 1.0 | 207.31 | 0.279 | 201.16 | 3.237 | 201.50 | 3.074 | 208.53 | 0.308 | 207.89 |
| | 2.0 | 221.00 | 0.629 | 213.79 | 3.871 | 214.47 | 3.566 | 222.41 | 0.004 | 222.40 |
| | 5.0 | 230.69 | 1.305 | 222.82 | 4.673 | 224.14 | 4.107 | 232.68 | 0.453 | 233.74 |
| | 10.0 | 232.98 | 1.617 | 225.06 | 4.964 | 226.75 | 4.248 | 235.54 | 0.536 | 236.81 |
| ARE (%) | | | 0.438 | | 1.814 | | 1.661 | | 0.189 | |

**Table 7** Comparison of the mid-span vertical displacement solutions [$10^{-2}$ mm] (CC).

| Type | *p* | DTS [23] | | DQ3D [30] | | MQ3D [30] | | FEQ3D | | Q4 |
|---|---|---|---|---|---|---|---|---|---|---|
| | | Disp | RE (%) | Disp | RE (%) | Disp | RE (%) | Disp | RE (%) | |
| A | 0.0 | 18.454 | 0.556 | 18.259 | 0.505 | 18.403 | 0.278 | 18.353 | 0.005 | 18.352 |
| | 0.5 | 27.109 | 0.702 | 26.825 | 0.353 | 27.019 | 0.368 | 26.926 | 0.022 | 26.920 |
| | 1.0 | 34.321 | 0.654 | 33.962 | 0.400 | 34.205 | 0.314 | 34.064 | 0.100 | 34.098 |
| | 2.0 | 44.649 | 0.593 | 44.173 | 0.479 | 44.524 | 0.311 | 44.290 | 0.216 | 44.386 |
| | 5.0 | 58.444 | 0.288 | 57.796 | 0.823 | 58.466 | 0.326 | 58.003 | 0.468 | 58.276 |
| | 10.0 | 65.241 | 0.092 | 64.509 | 1.214 | 65.516 | 0.329 | 64.824 | 0.730 | 65.301 |
| B | 0.0 | 35.150 | 0.077 | 34.283 | 2.541 | 34.764 | 1.174 | 35.320 | 0.407 | 35.177 |
| | 0.5 | 40.825 | 0.130 | 39.564 | 2.963 | 40.013 | 1.862 | 41.136 | 0.893 | 40.772 |
| | 1.0 | 44.081 | 0.120 | 42.587 | 3.504 | 43.155 | 2.218 | 44.578 | 1.006 | 44.134 |
| | 2.0 | 47.284 | 0.724 | 45.571 | 4.321 | 46.448 | 2.480 | 47.991 | 0.760 | 47.629 |
| | 5.0 | 49.929 | 1.755 | 48.071 | 5.411 | 49.503 | 2.593 | 50.880 | 0.116 | 50.821 |
| | 10.0 | 50.816 | 2.352 | 48.940 | 5.956 | 50.699 | 2.577 | 51.940 | 0.193 | 52.040 |
| ARE (%) | | | 0.670 | | 2.373 | | 1.236 | | 0.410 | |

Regarding the distributions of transverse normal stress, as shown in **Fig. 6** and **Fig. 9**, the results from DQ3D and MQ3D exhibit significant deviations from the Q4 reference solution. In contrast, the FEQ3D's results demonstrate excellent agreement with the Q4 solution. These findings not only validate the correctness and computational advantages of the proposed element formulation but also highlight the critical influence of employing equilibrium-based transverse normal stress expression on the accuracy of stress predictions.

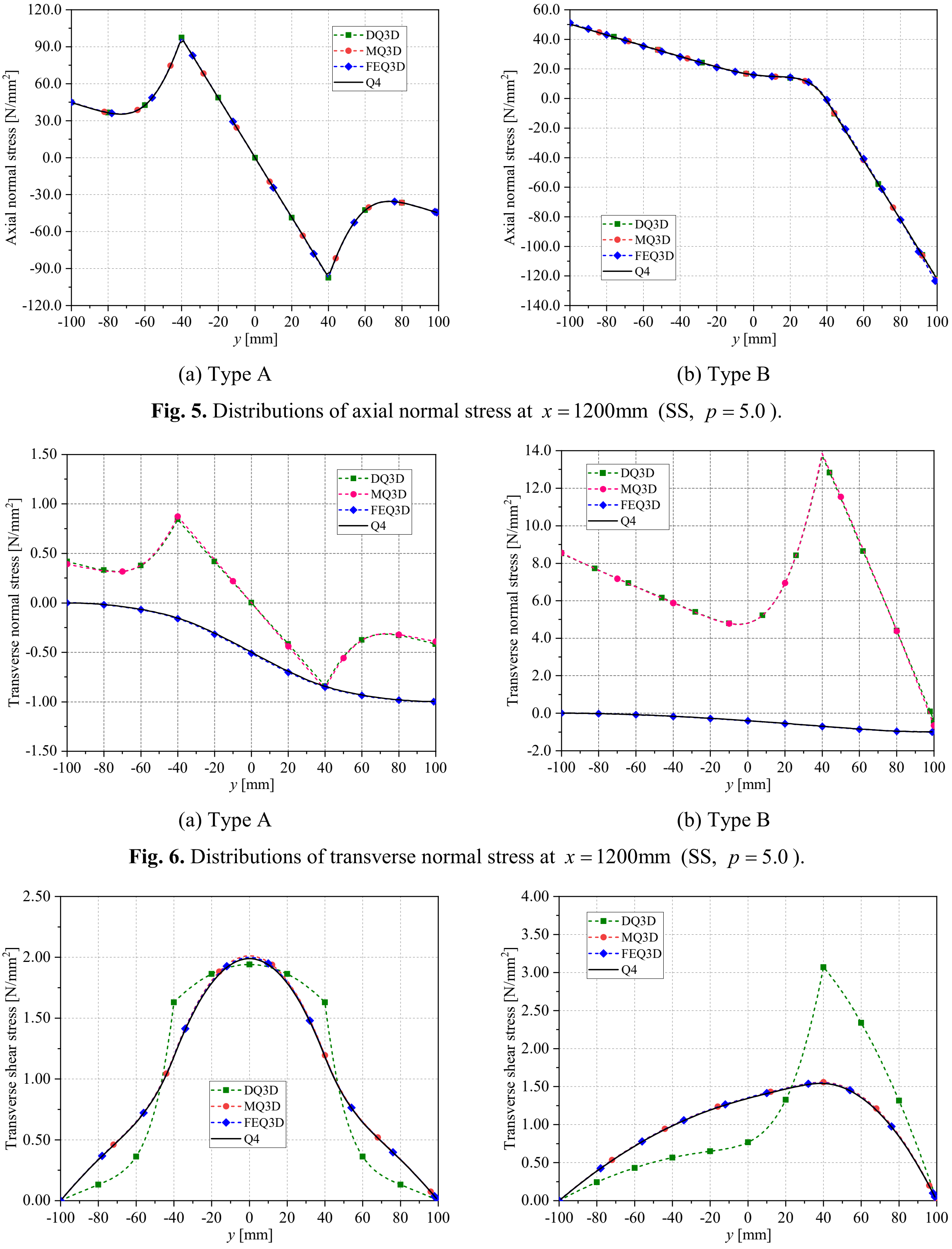

(a) Type A (b) Type B

**Fig. 5.** Distributions of axial normal stress at $x = 1200\text{mm}$ (SS, $p = 5.0$ ).

(a) Type A (b) Type B

**Fig. 6.** Distributions of transverse normal stress at $x = 1200\text{mm}$ (SS, $p = 5.0$ ).

(a) Type A (b) Type B

**Fig. 7.** Distributions of transverse shear stress at $x = 1200\text{mm}$ (SS, $p = 5.0$ ).

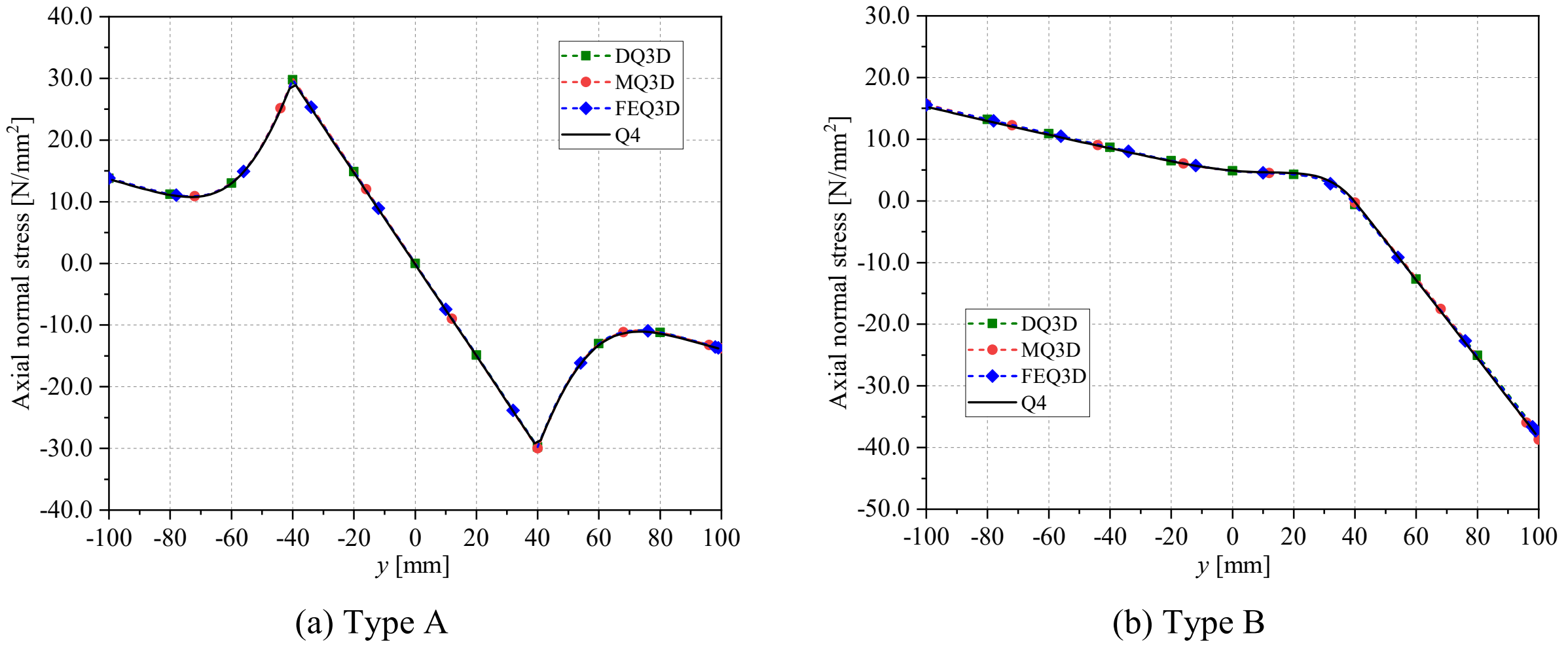

(a) Type A (b) Type B

**Fig. 8.** Distributions of axial normal stress at $x = 1200\text{mm}$ (CC, $p = 5.0$ ).

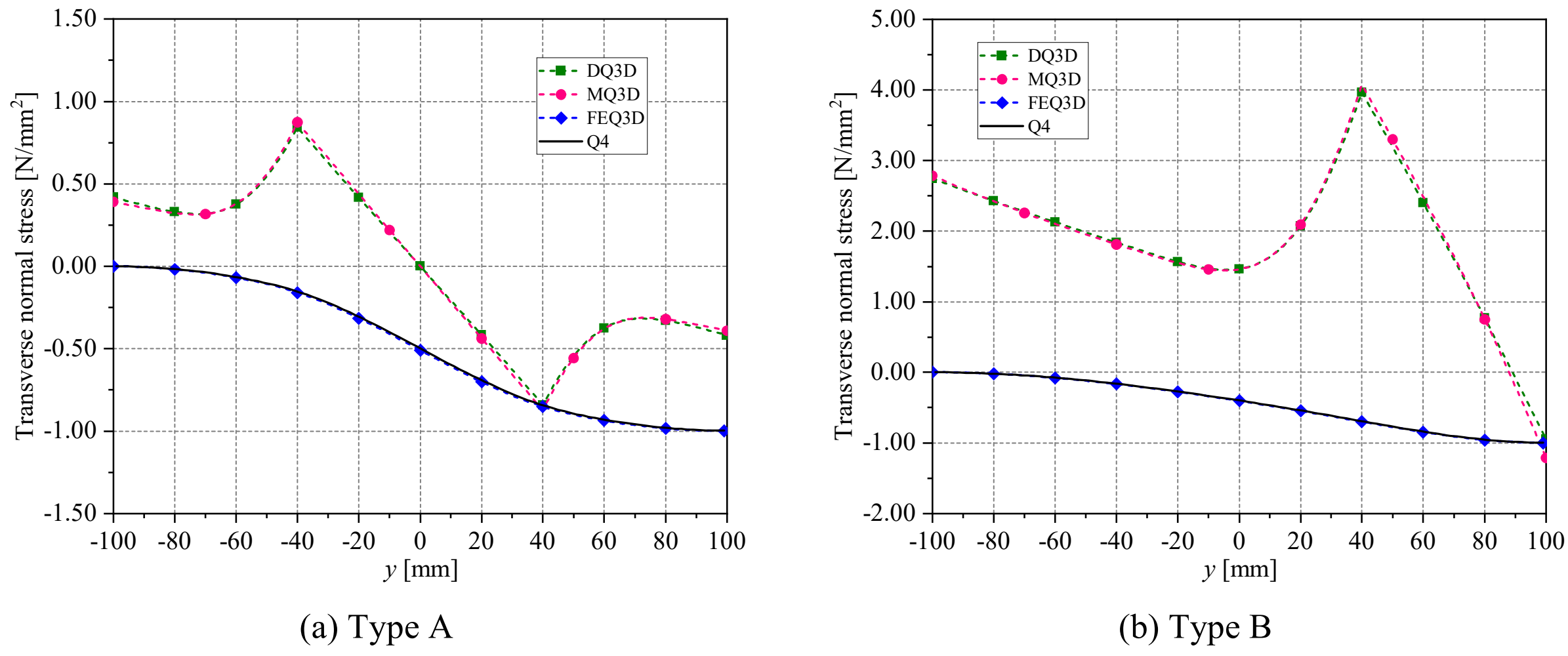

(a) Type A (b) Type B

**Fig. 9.** Distributions of transverse normal stress at $x = 1200\text{mm}$ (CC, $p = 5.0$ ).

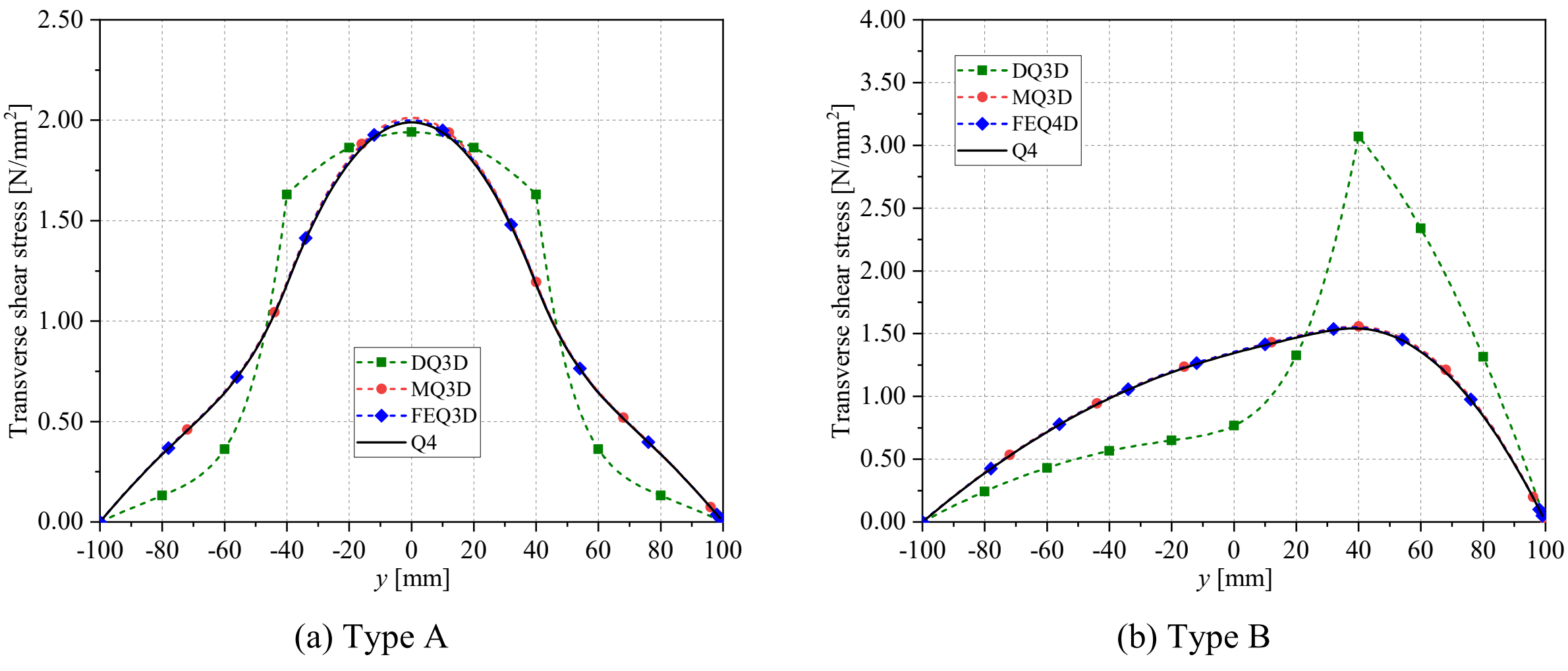

(a) Type A (b) Type B

**Fig. 10.** Distributions of transverse shear stress at $x = 1200\text{mm}$ (CC, $p = 5.0$ ).

For quasi-3D beam models, the inability to accurately predict transverse normal stress fundamentally compromises the determination of the actual stress state within the material. This deficiency significantly diminishes the practical utility of quasi-3D beam formulation. From this perspective, the proposed quasi-3D exact finite element demonstrates remarkable capability in simultaneously achieving accurate predictions for all three stress components, thereby substantially enhancing the applicability and engineering relevance of quasi-3D beam models.

### *5.3. Impact of the higher-order functions*

Under the framework of DQ3D and FEQ3D, this section investigates the influence of the higher-order function $f(y)$ on solution accuracy. The three forms of higher-order function listed in **Table 1** are considered, and the Third-order polynomial function, Sinusoidal function, and Exponential function are denoted by "-T," "-S," and "-E," respectively. With the use of the three higher-order functions, the displacement results obtained using DQ3D are presented in **Table 8** and **Table 9**, while those obtained using FEQ3D are shown in **Table 10** and **Table 11**.

The results in **Table 8** and **Table 9** indicate that, within the conventional displacement-based quasi-3D beam element framework, the choice of higher-order function has a notable impact on the element accuracy, particularly for the case with the Type A material model. Among the three higher-order functions considered, the element model employing a third-order polynomial form (DQ3D-T) achieves the highest accuracy. For the Type B material model, the average relative error (ARE) of DQ3D is approximately 4%, and no significant improvement can be achieved by adopting different higher-order functions. For the FEQ3D framework, as shown in **Table 10** and **Table 11**, the influence of the higher-order function form on element accuracy is relatively minor. According to the results, for symmetric material distributions (Type A), the element model with the third-order polynomial higher-order function (FEQ3D-T) exhibits the highest accuracy. For asymmetric material distributions (Type B), the differences in accuracy among the elements with the three higher-order functions are minimal, and the sinusoidal higher-order function yields the highest average accuracy. Overall, regardless of the higher-order function used, FEQ3D demonstrates superior solution accuracy compared to DQ3D.

For both frameworks of DQ3D and FEQ3D, under the clamped-clamped boundary condition, the stress distributions through the cross-section located at $x = 1200\text{mm}$, obtained using element models based on the three higher-order functions, are presented in **Fig. 11**-**Fig. 13**, corresponding to the distributions of axial normal stress, transverse normal stress, and transverse shear stress, respectively. It can be observed that for the axial normal stress distribution, the influence of the higher-order function is negligible. However, for the transverse shear stress and transverse normal stress distributions, the impact of higher-order function form varies depending on the element framework. For the proposed force-based quasi-3D beam exact finite element model, the influence of the higher-order function is minimal, whether for transverse shear stress or transverse normal stress. In contrast, for the conventional displacement-based quasi-3D beam element model, the higher-order function form has a noticeable effect on the transverse shear stress distribution and a significant impact on the transverse normal stress distribution. Notably, due to the substantial differences in the second-order derivatives of the various higher-order functions, the transverse normal stress distribution obtained using DQ3D varies significantly with the higher-order function form. Overall, despite employing different higher-order function forms, the accuracy of stress solutions obtained using DQ3D cannot be effectively improved. In contrast, under the framework of FEQ3D, relatively accurate stress solutions can be obtained regardless of the higher-order function employed.

**Table 8** Comparison of the mid-span vertical displacements for DQ3D with different $f(y)$ [$10^{-2}$ mm] (Type A).

| Boundary | $p$ | DQ3D-T | | DQ3D-S | | DQ3D-E | | Q4 |
|---|---|---|---|---|---|---|---|---|
| | | Disp | RE (%) | Disp | RE (%) | Disp | RE (%) | |
| SS | 0.0 | 84.110 | 0.202 | 83.995 | 0.339 | 83.652 | 0.745 | 84.280 |
| | 0.5 | 126.29 | 0.052 | 126.15 | 0.170 | 125.70 | 0.522 | 126.36 |
| | 1.0 | 161.61 | 0.013 | 161.43 | 0.124 | 160.90 | 0.454 | 161.63 |
| | 2.0 | 212.37 | 0.048 | 212.13 | 0.068 | 211.41 | 0.404 | 212.27 |
| | 5.0 | 280.41 | 0.047 | 280.00 | 0.101 | 278.87 | 0.502 | 280.28 |
| | 10.0 | 313.96 | 0.029 | 313.43 | 0.197 | 312.02 | 0.645 | 314.05 |
| CC | 0.0 | 18.259 | 0.505 | 18.223 | 0.701 | 18.144 | 1.136 | 18.352 |
| | 0.5 | 26.825 | 0.353 | 26.782 | 0.513 | 26.684 | 0.878 | 26.920 |
| | 1.0 | 33.962 | 0.400 | 33.911 | 0.547 | 33.796 | 0.885 | 34.098 |
| | 2.0 | 44.173 | 0.479 | 44.106 | 0.631 | 43.955 | 0.971 | 44.386 |
| | 5.0 | 57.796 | 0.823 | 57.684 | 1.015 | 57.449 | 1.419 | 58.276 |
| | 10.0 | 64.509 | 1.214 | 64.364 | 1.436 | 64.067 | 1.890 | 65.301 |
| ARE (%) | | | 0.347 | | 0.487 | | 0.871 | |

**Table 9** Comparison of the mid-span vertical displacements for DQ3D with different $f(y)$ [$10^{-2}$ mm] (Type B).

| Boundary | $p$ | DQ3D-T | | DQ3D-S | | DQ3D-E | | Q4 |
|---|---|---|---|---|---|---|---|---|
| | | Disp | RE (%) | Disp | RE (%) | Disp | RE (%) | |
| SS | 0.0 | 163.27 | 1.928 | 163.06 | 2.057 | 162.47 | 2.410 | 166.48 |
| | 0.5 | 187.63 | 2.705 | 187.16 | 2.946 | 186.33 | 3.377 | 192.84 |
| | 1.0 | 201.16 | 3.237 | 200.56 | 3.525 | 199.60 | 3.987 | 207.89 |
| | 2.0 | 213.79 | 3.871 | 213.12 | 4.171 | 212.10 | 4.632 | 222.40 |
| | 5.0 | 222.82 | 4.673 | 222.25 | 4.918 | 221.33 | 5.309 | 233.74 |
| | 10.0 | 225.06 | 4.964 | 224.60 | 5.155 | 223.82 | 5.486 | 236.81 |
| CC | 0.0 | 34.283 | 2.541 | 34.206 | 2.761 | 34.061 | 3.174 | 35.177 |
| | 0.5 | 39.564 | 2.963 | 39.440 | 3.267 | 39.247 | 3.740 | 40.772 |
| | 1.0 | 45.571 | 3.504 | 42.446 | 3.824 | 42.235 | 4.302 | 44.134 |
| | 2.0 | 45.522 | 4.321 | 45.438 | 4.600 | 45.236 | 5.024 | 47.629 |
| | 5.0 | 48.071 | 5.411 | 48.007 | 5.537 | 47.878 | 5.790 | 50.821 |
| | 10.0 | 48.940 | 5.956 | 48.934 | 5.968 | 48.868 | 6.096 | 52.040 |
| ARE (%) | | | 3.839 | | 4.061 | | 4.444 | |

**Table 10** Comparison of the mid-span vertical displacements for FEQ3D with different $f(y)$ [$10^{-2}$ mm] (Type A).

| Boundary | *p* | FEQ3D-T | | FEQ3D-S | | FEQ3D-E | | Q4 |
|---|---|---|---|---|---|---|---|---|
| | | Disp | RE (%) | Disp | RE (%) | Disp | RE (%) | |
| SS | 0.0 | 84.109 | 0.203 | 84.118 | 0.193 | 84.098 | 0.217 | 84.280 |
| | 0.5 | 126.30 | 0.047 | 126.30 | 0.047 | 126.24 | 0.096 | 126.36 |
| | 1.0 | 161.61 | 0.012 | 161.61 | 0.015 | 161.51 | 0.074 | 161.63 |
| | 2.0 | 212.39 | 0.057 | 212.38 | 0.050 | 212.24 | 0.015 | 212.27 |
| | 5.0 | 280.48 | 0.071 | 280.44 | 0.057 | 280.24 | 0.013 | 280.28 |
| | 10.0 | 314.07 | 0.006 | 314.02 | 0.011 | 313.79 | 0.084 | 314.05 |
| CC | 0.0 | 18.353 | 0.005 | 18.344 | 0.043 | 18.326 | 0.142 | 18.352 |
| | 0.5 | 26.926 | 0.022 | 26.913 | 0.024 | 26.887 | 0.124 | 26.920 |
| | 1.0 | 34.064 | 0.100 | 34.048 | 0.147 | 34.013 | 0.248 | 34.098 |
| | 2.0 | 44.290 | 0.216 | 44.264 | 0.274 | 44.216 | 0.384 | 44.386 |
| | 5.0 | 58.003 | 0.468 | 57.953 | 0.554 | 57.875 | 0.688 | 58.276 |
| | 10.0 | 64.824 | 0.730 | 64.755 | 0.836 | 64.657 | 0.987 | 65.301 |
| ARE (%) | | | 0.161 | | 0.187 | | 0.256 | |

**Table 11** Comparison of the mid-span vertical displacements for FEQ3D with different $f(y)$ [$10^{-2}$ mm] (Type B).

| Boundary | *p* | FEQ3D-T | | FEQ3D-S | | FEQ3D-E | | Q4 |
|---|---|---|---|---|---|---|---|---|
| | | Disp | RE (%) | Disp | RE (%) | Disp | RE (%) | |
| SS | 0.0 | 166.80 | 0.192 | 166.83 | 0.210 | 166.76 | 0.171 | 166.48 |
| | 0.5 | 193.57 | 0.379 | 193.52 | 0.351 | 193.41 | 0.295 | 192.84 |
| | 1.0 | 208.53 | 0.308 | 208.41 | 0.248 | 208.25 | 0.175 | 207.89 |
| | 2.0 | 222.41 | 0.004 | 222.26 | 0.061 | 222.14 | 0.118 | 222.40 |
| | 5.0 | 232.68 | 0.453 | 232.73 | 0.433 | 232.80 | 0.402 | 233.74 |
| | 10.0 | 235.54 | 0.536 | 235.72 | 0.459 | 235.92 | 0.374 | 236.81 |
| CC | 0.0 | 35.320 | 0.407 | 35.308 | 0.373 | 35.278 | 0.288 | 35.177 |
| | 0.5 | 41.136 | 0.893 | 41.116 | 0.844 | 41.082 | 0.759 | 40.772 |
| | 1.0 | 44.578 | 1.006 | 44.543 | 0.926 | 44.504 | 0.838 | 44.134 |
| | 2.0 | 47.991 | 0.760 | 47.972 | 0.719 | 47.961 | 0.696 | 47.629 |
| | 5.0 | 50.880 | 0.116 | 50.973 | 0.299 | 51.074 | 0.498 | 50.821 |
| | 10.0 | 51.940 | 0.193 | 52.115 | 0.144 | 52.291 | 0.482 | 52.040 |
| ARE (%) | | | 0.437 | | 0.422 | | 0.425 | |

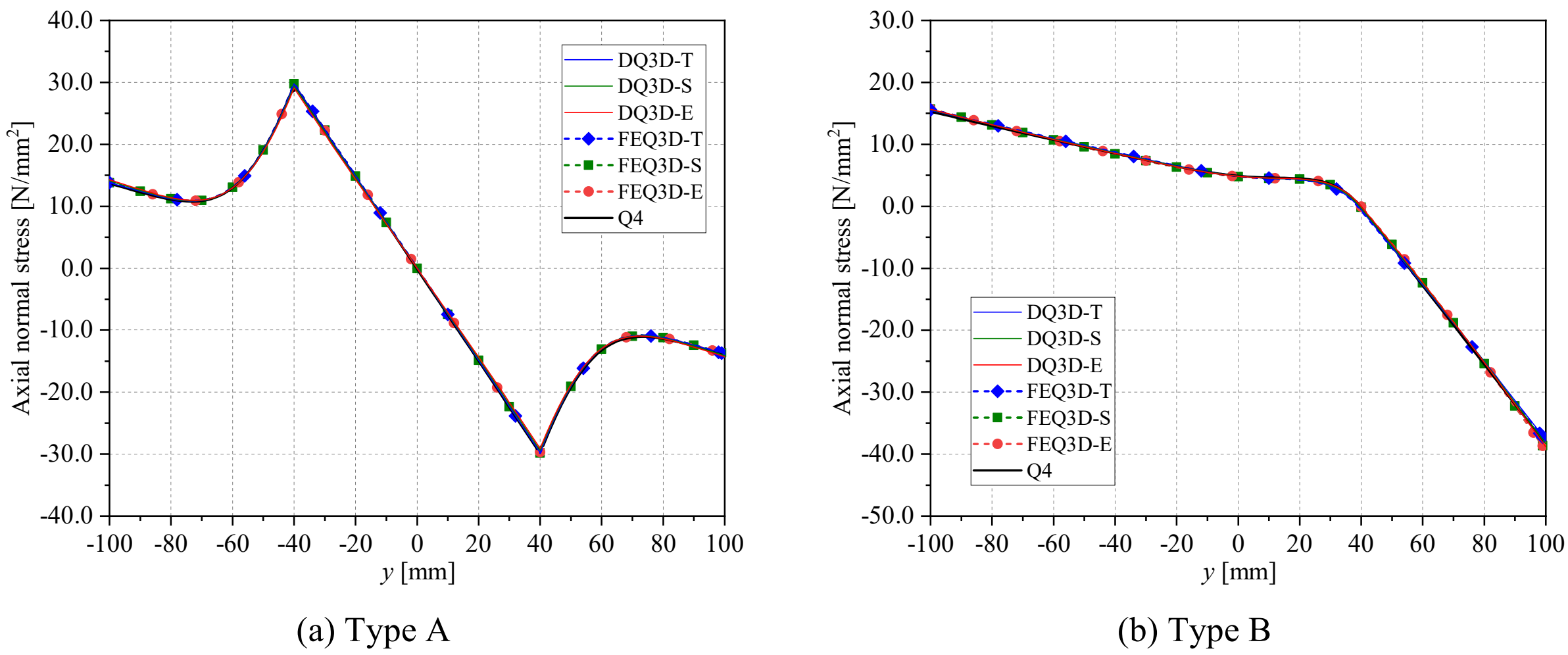

(a) Type A (b) Type B

**Fig. 11.** Distributions of axial normal stress for different $f(y)$ (CC, $p = 5.0$, $x = 1200\text{mm}$).

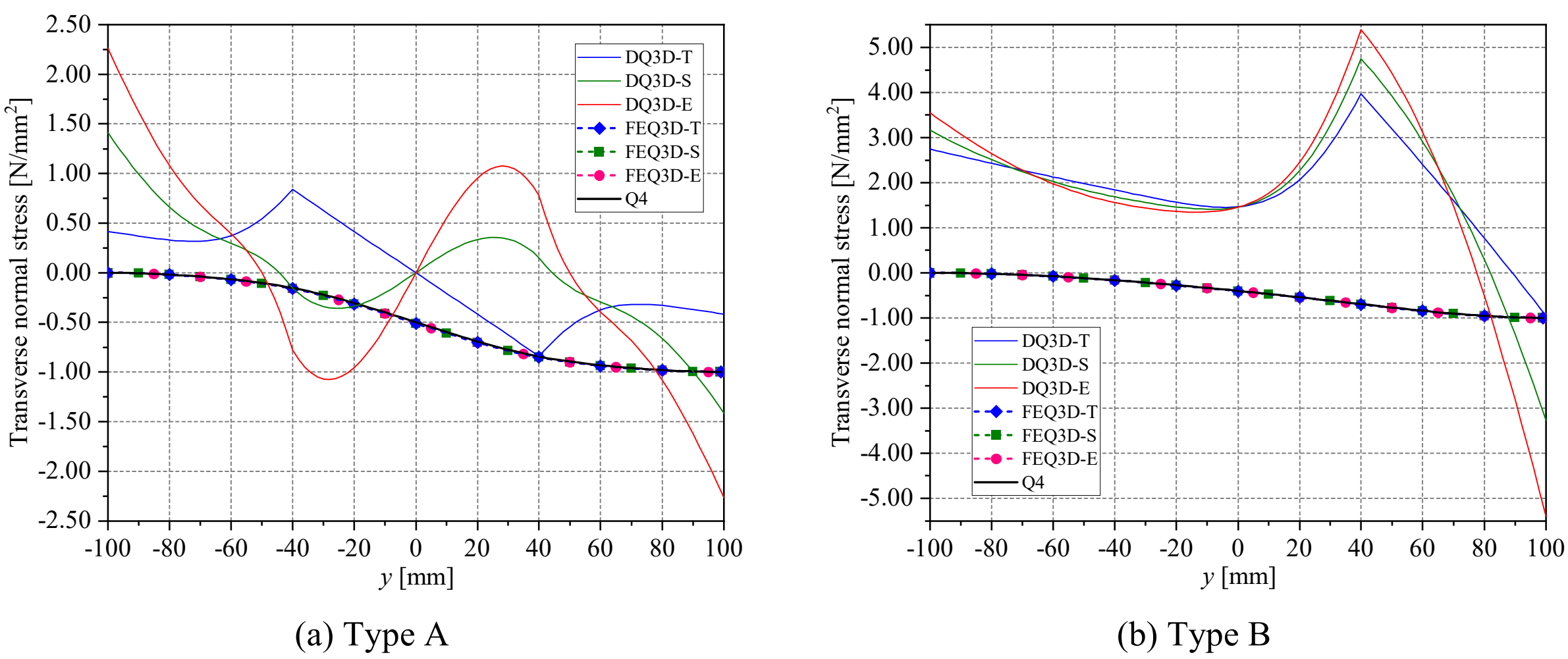

(a) Type A (b) Type B

**Fig. 12.** Distributions of transverse normal stress for different $f(y)$ (CC, $p = 5.0$, $x = 1200\text{mm}$).

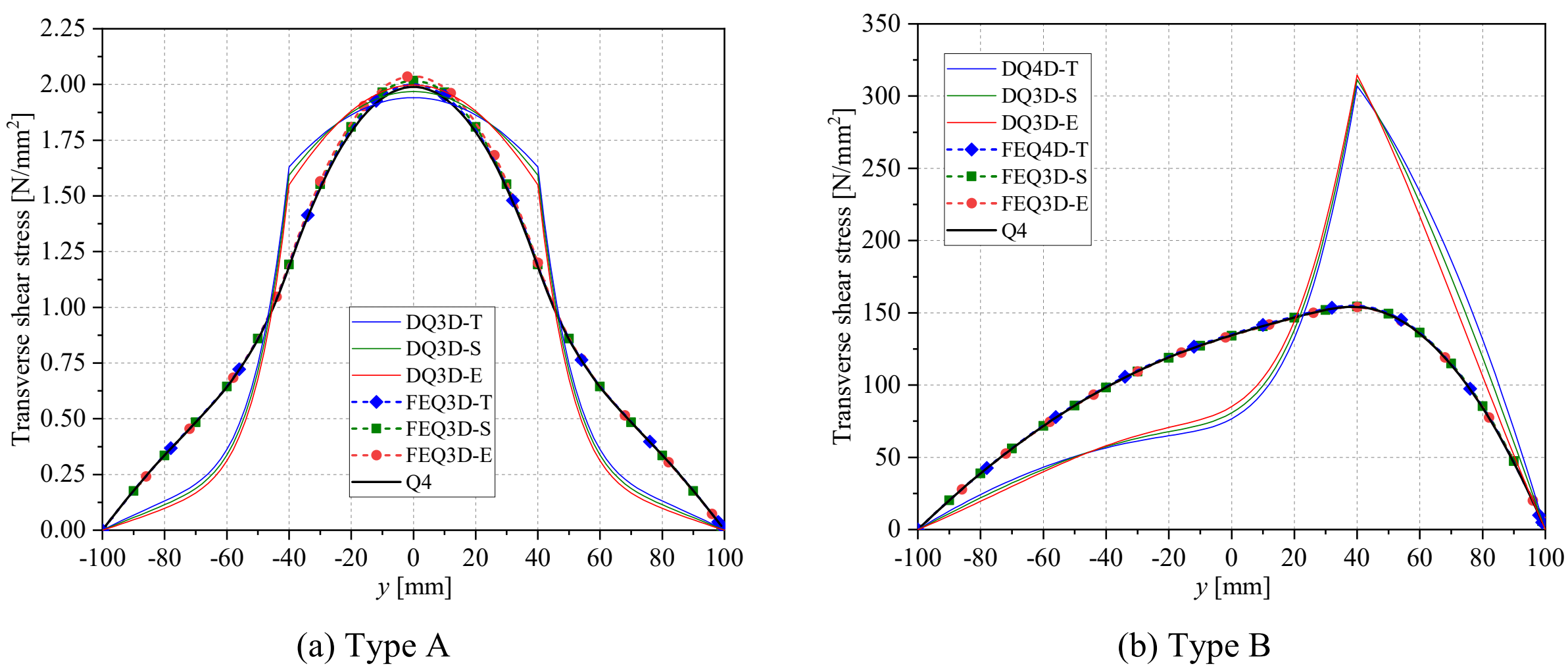

(a) Type A (b) Type B

**Fig. 13.** Distributions of transverse shear stress for different $f(y)$ (CC, $p = 5.0$, $x = 1200\text{mm}$).

*5.4. Impact of the distributed loads*

This section further investigates the performance of FEQ3D (with third-order polynomial higher-order function) under gradient-distributed loads. Firstly, the linearly distributed load is considered, where $q_2 = 0$. In this investigation, the mean intensity of distribution load is set to $\overline{q} = -50\,\mathrm{N/mm}$, then the load intensity at the starting end can be determined according to the setting of $q_1$ as

$$q_0 = \overline{q} - \frac{q_1 L}{2} \tag{120}$$

In the numerical validation, two computational approaches are considered: (1) The stiffness equations of the structure are assembled from the element stiffness equations with distributed load included, which has been presented in Eq. (90); and (2) The stiffness equations of the structure are assembled from the element stiffness matrix derived from FEQ3D without distributed load, while the equivalent nodal loads are introduced to form the nodal load vector. For the second approach, as shown in **Fig. 14**, the equivalent nodal load for the *i*-th internal node is determined as follows for the mesh of *n* elements

$$P_i = q\left(\frac{L}{n} i\right)\frac{L}{n} \tag{121}$$

For the sake of distinction, the two approaches mentioned above are denoted as FEQ3D and FEQ3D$^{\mathrm{C}}$, respectively.

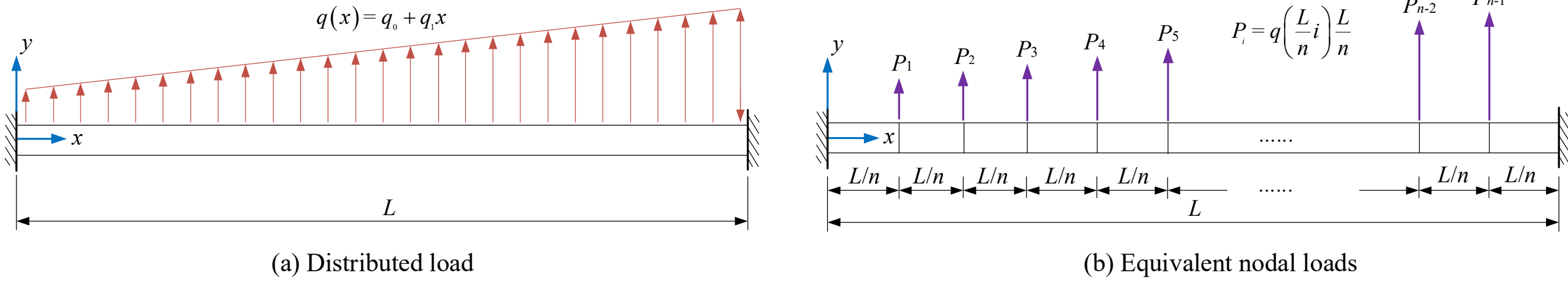

(a) Distributed load (b) Equivalent nodal loads

**Fig. 14.** Setting of equivalent nodal loads.

The study examines the maximum vertical displacement of beams under four predefined gradient configurations, including $q_1 = -0.02\,\mathrm{N/mm^2}$, $q_1 = -0.05\,\mathrm{N/mm^2}$, $q_1 = -0.10\,\mathrm{N/mm^2}$ and $q_1 = -0.20\,\mathrm{N/mm^2}$, with results presented in **Table 12** and **Table 13** for Type A and Type B material models, respectively, considering various mesh refinement levels under two computational approaches. The findings demonstrate that for gradient-distributed loads, the proposed formulation achieves converged displacement results with just a single element, eliminating the need for mesh refinement. When employing the second computational approach (equivalent nodal loading), convergence can be achieved through progressive mesh refinement, and the final converged values match those obtained from the first computational approach. The results confirms the accuracy and robustness of the present formulation in handling gradient-distributed loads. Meanwhile, the investigation of the vertical displacement distribution, as shown in **Fig. 15**, reveals that the maximum displacement shifts toward the side with higher load intensity rather than occurring at the mid-span, which is physically consistent with expectations due to the non-uniform load distribution.

**Table 12** Maximum displacements from FEQ3D under gradient-distributed loads [$10^{-2}$ mm] (CC, Type A, $p = 5.0$ ).

| Number of elements | $q_1 = -0.02$ | | $q_1 = -0.05$ | | $q_1 = -0.10$ | | $q_1 = -0.20$ | |
|---|---|---|---|---|---|---|---|---|
| | FEQ3D$^C$ | FEQ3D | FEQ3D$^C$ | FEQ3D | FEQ3D$^C$ | FEQ3D | FEQ3D$^C$ | FEQ3D |
| 1 | - | 58.055 | - | 58.325 | - | 59.225 | - | 62.590 |
| 2 | 57.674 | 58.055 | 57.674 | 58.325 | 57.674 | 59.225 | 57.674 | 62.590 |
| 4 | 57.817 | | 57.817 | | 57.817 | | 57.817 | |
| 8 | 57.916 | | 57.916 | | 57.916 | | 61.665 | |
| 16 | 57.978 | | 57.978 | | 59.155 | | 62.097 | |
| 32 | 58.001 | | 58.311 | | 59.184 | | 62.578 | |
| 64 | 58.044 | | 58.313 | | 59.248 | | 62.582 | |
| 128 | 58.051 | | 58.321 | | 59.248 | | 62.582 | |
| 256 | 58.054 | | 58.324 | | 59.255 | | 62.589 | |
| 512 | 58.055 | | 58.325 | | 59.255 | | 62.590 | |
| 1024 | 58.055 | | 58.325 | | 59.255 | | 62.590 | |
| Converged | 58.055 | 58.055 | 58.325 | 58.325 | 59.255 | 59.255 | 62.590 | 62.590 |

**Table 13** Maximum displacements from FEQ3D under gradient-distributed loads [$10^{-2}$ mm] (CC, Type B, $p = 5.0$ ).

| Number of elements | $q_1 = -0.02$ | | $q_1 = -0.05$ | | $q_1 = -0.10$ | | $q_1 = -0.20$ | |
|---|---|---|---|---|---|---|---|---|
| | FEQ3D$^C$ | FEQ3D | FEQ3D$^C$ | FEQ3D | FEQ3D$^C$ | FEQ3D | FEQ3D$^C$ | FEQ3D |
| 1 | - | 50.932 | - | 51.203 | - | 52.136 | - | 55.442 |
| 2 | 50.724 | 50.932 | 50.724 | 51.203 | 50.724 | 52.136 | 50.724 | 55.442 |
| 4 | 50.808 | | 50.808 | | 50.808 | | 50.808 | |
| 8 | 50.851 | | 50.851 | | 50.851 | | 54.932 | |
| 16 | 50.873 | | 50.873 | | 52.102 | | 55.003 | |
| 32 | 50.879 | | 51.200 | | 52.111 | | 55.436 | |
| 64 | 50.926 | | 51.200 | | 52.112 | | 55.437 | |
| 128 | 50.926 | | 51.200 | | 52.136 | | 55.437 | |
| 256 | 50.932 | | 51.202 | | 52.136 | | 55.442 | |
| 512 | 50.932 | | 51.203 | | 52.136 | | 55.442 | |
| 1024 | 50.932 | | 51.203 | | 52.136 | | 55.442 | |
| Converged | 50.932 | 50.932 | 51.203 | 51.203 | 52.136 | 52.136 | 55.442 | 55.442 |

The internal force distributions under varying gradient conditions (using CC boundary as an example) are presented in **Fig. 16**. The analysis shows that the shapes of the two bending moment components remain similar (with varying magnitudes). The moment distribution curves exhibit higher-order curvature (beyond second-order) that becomes increasingly complex with higher load gradients. Regarding the shear force, the linear distribution assumption is invalidated within the main span ( $x \in [100\text{mm}, 1900\text{mm}]$ , the shear force distribution is complex near the boundaries due to the boundary effects). Moreover, it should be noted that the stress resultant corresponding to transverse stretching deformation, $R(x)$, varies linearly (within the main span), directly correlating with the load intensity. Generally, higher load intensity produces proportionally higher values of transverse stretching stress resultant. As a supplement, the stress resultants corresponding to transverse stretching deformation obtained by FEQ3D$^C$ under various mesh refinements are presented in **Fig. 17**. The results demonstrate that with coarse meshes, the stress resultant distributions exhibit oscillatory behavior. As mesh refinement increases, the solution converges toward physically reasonable results.

Collectively, these findings confirm that the proposed formulation effectively captures gradient load effects, enabling exact finite element modeling of beams subjected to non-uniform load distributions. The consistency between results from both computational approaches (with sufficient mesh refinement) further validates the reliability of the proposed approach for engineering applications involving gradient-distributed loads.

In addition, the solution accuracy of FEQ3D under a quadratic distributed load is further investigated. For the two settings of load parameters (Setting 1: $q_0 = 0$, $q_1 = -0.01\,\mathrm{N/mm^2}$, $q_2 = 5\times10^{-5}\,\mathrm{N/mm^3}$; Setting 2: $q_0 = -50\,\mathrm{N/mm}$, $q_1 = -0.02\,\mathrm{N/mm^2}$, $q_2 = 8\times10^{-5}\,\mathrm{N/mm^3}$), the maximum vertical displacement of beams (CC, Type B) obtained by FEQ3D and FEQ3D$^C$ are presented in **Table 14**. The results demonstrate that FEQ3D achieves accurate solutions for the quadratic distributed loads by using only one beam element. In fact, the proposed method is not limited to cases involving polynomial-distributed loads. Still, it is also applicable to more complex distributed loading conditions, provided that the corresponding formulation can be updated based on the characteristics of the load distribution.

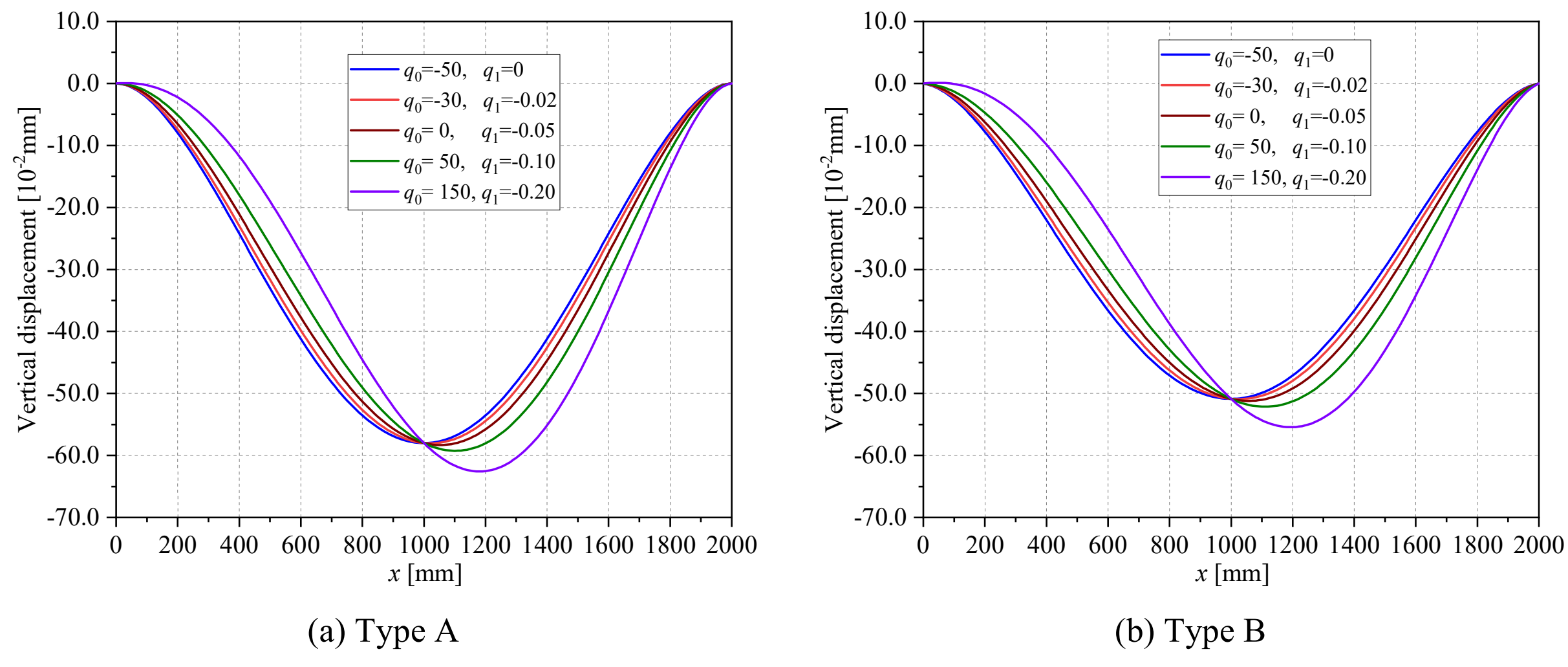

(a) Type A (b) Type B

**Fig. 15.** Distributions of vertical displacement along $x$ coordinate (FEQ3D, CC, $p = 5.0$).

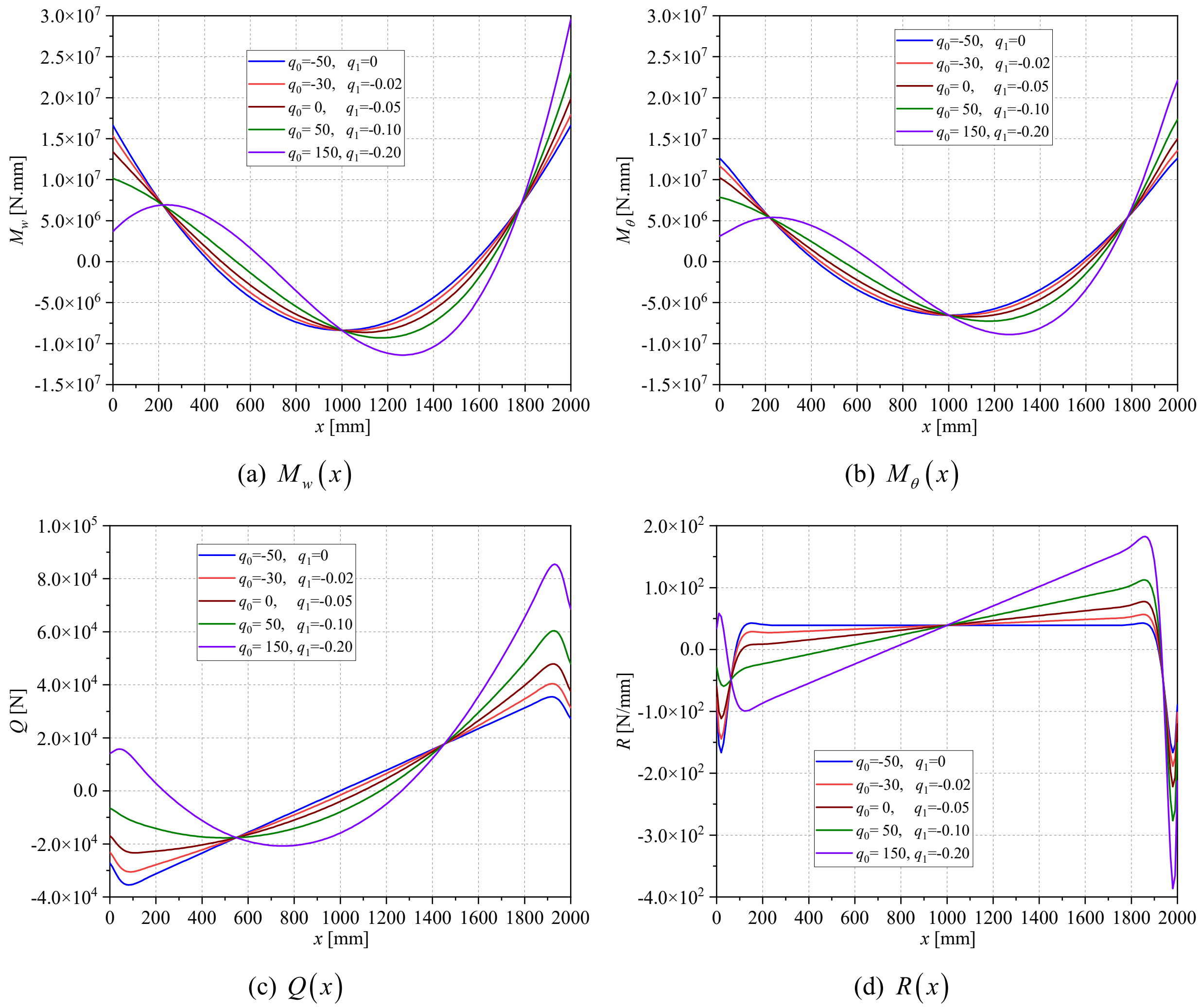

(a) $M_w(x)$ (b) $M_\theta(x)$

(c) $Q(x)$ (d) $R(x)$

**Fig. 16.** Distributions of internal forces along $x$ coordinate (FEQ3D, CC, Type B, $p = 5.0$ ).

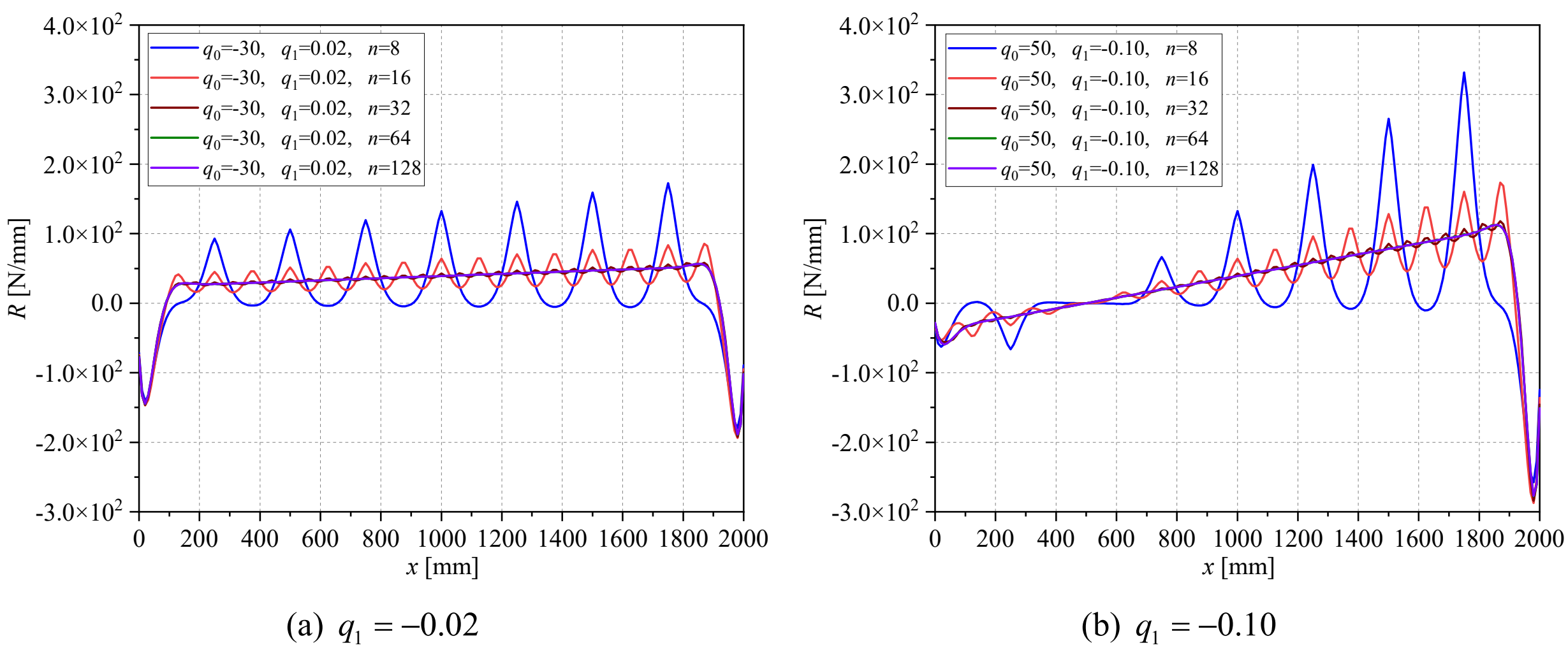

(a) $q_1 = -0.02$ (b) $q_1 = -0.10$

**Fig. 17.** Distributions of $R(x)$ from FEQ3D[C] (CC, Type B, $p = 5.0$ ).

**Table 14** Maximum displacements from FEQ3D under quadratic-distributed loads [$10^{-2}$ mm] (CC, Type B, $p = 5.0$ ).

| Number of elements | Setting 1: $q_0 = 0, q_1 = -0.01, q_2 = 5\times10^{-5}$ | | Setting 2: $q_0 = -50, q_1 = -0.02, q_2 = 8\times10^{-5}$ | |
|---|---|---|---|---|
| | FEQ3D$^{C}$ | FEQ3D | FEQ3D$^{C}$ | FEQ3D |
| 1 | - | 43.929 | - | 161.93 |
| 2 | 50.724 | 43.929 | 172.46 | 161.93 |
| 8 | 43.969 | | 161.86 | |
| 32 | 43.929 | | 161.87 | |
| 128 | 43.929 | | 161.93 | |
| 512 | 43.929 | | 161.93 | |
| Converged | 43.929 | 43.929 | 161.93 | 161.93 |

## 6. Conclusions

This study presents an exact finite element formulation of quasi-3D beam based on analytical solutions of internal force fields, aiming to achieve high-fidelity analysis of functionally graded sandwich beams. Unlike conventional displacement-based elements that rely on approximate interpolation functions or existing exact finite element methods requiring closed-form solutions of generalized displacements, the novelty of the present work lies in constructing an exact beam element from analytical expressions of internal forces, integrating the force-based beam formulation with the exact finite element framework to address the limitations of existing approaches. The main conclusions of this work are summarized as follows:

(1) The proposed approach of constructing exact finite elements based on analytical internal force fields is validated as correct. The proposed formulation establishes an efficient and high-precision finite element model that effectively eliminates discretization errors and yields a symmetric element stiffness matrix.

(2) The proposed exact finite element formulation successfully incorporates quasi-3D beam theory, enabling the development of a high-accuracy beam element model capable of capturing transverse stretching deformation, which is considered as a critical feature for accurately modeling FG sandwich beams.

(3) By introducing the modified cross-sectional stiffness matrices that accounts for equilibrium-based stress distributions, the proposed quasi-3D beam exact finite element achieves simultaneous accuracy in both displacement and stress solutions for FG sandwich beams under static loading conditions.

(4) Within the proposed exact finite element framework of quasi-3D beam, the element accuracy exhibits minimal sensitivity to the choice of higher-order functions used to describe the distribution of out-of-plane axial displacement through the thickness. This characteristic demonstrates superior stability compared to conventional displacement-based elements, where the accuracy is more sensitive to the selection of higher-order functions.

(5) The proposed formulation accurately considers the influence of gradient-distributed loads. For more complex distributed loads, corresponding element stiffness equations can be systematically derived following the same theoretical framework, ensuring broad applicability for engineering problems involving non-uniform loading.

Despite the significant advantages of the proposed element, it is also crucial to acknowledge its limitations. First, the formulation of the present exact element is predicated on the availability of analytical expressions for the internal force fields. For more complex scenarios, such as beams with material properties varying along the *x*-axis, obtaining such analytical solutions becomes exceedingly difficult, which would pose a significant obstacle to the application of the current method. Furthermore, the extension of the current formulation to account for essential nonlinearities in practical engineering, including large deformations and material nonlinearity, constitutes a substantial challenge. Therefore, investigations addressing these limitations are planned as subsequent research directions.

**Acknowledgments**

The project is funded by the National Natural Science Foundation of China (Grant No. 52178209, Grant No. 51878299), Guangdong Basic and Applied Basic Research Foundation, China (Grant No. 2025A1515011664, Grant No. 2021A1515012280, Grant No. 2020A1515010611), Science and Technology Innovation Program from Water Resources of Guangdong Province (Grant No. 2025-02) and Innovation Team Project for Ordinary Universities in Guangdong Province (Grant No. 2023KCXTD005).

**Appendix. Details on the solution for the SS beam with the Type B material model using FEQ3D**

The following unit convention is adopted throughout this appendix: [Force: N; Length: mm]. Units are not appended to each value in the equations for the sake of simplicity.

For the beam with Type B material model, the cross-sectional stiffness matrices can be obtained according to Eq. (17) as

$$\mathbf{D}_{\mathrm{n}} = \begin{bmatrix} 2.0183\text{E}9 & 7.8027\text{E}10 & 6.3966\text{E}10 & -4.6816\text{E}6 \\ 7.8027\text{E}10 & 8.0860\text{E}12 & 6.4414\text{E}12 & -4.8516\text{E}8 \\ 6.3966\text{E}10 & 6.4414\text{E}12 & 5.1891\text{E}12 & -3.8648\text{E}8 \\ -4.6816\text{E}6 & -4.8516\text{E}8 & -3.8648\text{E}8 & 3.2344\text{E}5 \end{bmatrix} \tag{A1a}$$

$$D_{55} = 3.1308\text{E}8 \tag{A1b}$$

Further, according to the formulation given in **Sec. 4.2**, the modified cross-sectional stiffness matrices can be obtained as

$$\mathbf{D}_{\mathrm{n}}^{\mathrm{m}} = \begin{bmatrix} 1.9945\text{E}9 & 7.8216\text{E}10 & 6.3922\text{E}10 & -4.8075\text{E}6 \\ 7.8216\text{E}10 & 7.6540\text{E}12 & 6.0902\text{E}12 & -1.9712\text{E}8 \\ 6.3922\text{E}10 & 6.0902\text{E}12 & 4.9013\text{E}12 & -1.5235\text{E}8 \\ -4.8075\text{E}6 & -1.9712\text{E}8 & -1.5235\text{E}8 & 1.3142\text{E}5 \end{bmatrix} \tag{A2a}$$

$$D_{55}^{\mathrm{m}} = 2.7368\text{E}8 \tag{A2b}$$

Considering the distribution load $q(x) = -50$, the axial force $N(x)$ and bending moment $M_w(x)$ can be defined as

$$N(x) = c_0 \tag{A3}$$

$$M_w(x) = c_1 + c_2 x - 25x^2 \tag{A4}$$

Then, Eq. (27) can be rewritten as

$$M_{\theta,xxxx}(x) + b_1 M_{\theta,xx}(x) + b_0 M_\theta(x) = \frac{25 f_{23}}{f_{44}} x^2 - \frac{f_{23}}{f_{44}} c_2 x + \frac{50 f_{24}}{f_{44}} - \frac{f_{13}}{f_{44}} c_0 - \frac{f_{23}}{f_{44}} c_1 \tag{A5}$$

By solving the eigen-equation given in Eq. (29), two pairs of conjugate complex roots with opposite real parts and identical imaginary parts are generally obtained as

$$\lambda_{1,2} = 0.0307 \pm 0.0237i \tag{A6a}$$

$$\lambda_{3,4} = -0.0307 \pm 0.0237i \tag{A6b}$$

Then, the general solution of Eq. (A5) can be expressed as

$$\bar{M}_\theta\left(x\right) = c_3 e^{0.0307x}\cos\left(0.0237x\right) + c_4 e^{0.0307x}\sin\left(0.0237x\right) + c_5 e^{-0.0307x}\cos\left(0.0237x\right) + c_6 e^{-0.0307x}\sin\left(0.0237x\right) \tag{A7}$$

Further, by assuming the particular solution of Eq. (A5) as

$$M_\theta^*\left(x\right) = a_2 x^2 + a_1 x + a_0 \tag{A8}$$

The coefficients $a_0 \sim a_2$ can be determined by substituting Eq. (A8) into Eq. (A5), and then the internal force fields can be expressed as

$$\boldsymbol{\sigma}_{\mathrm{n}}\left(x\right) = \begin{bmatrix} 1 & 0 & 0 & 0 & 0 & 0 & 0 \\ 0 & 1 & x & 0 & 0 & 0 & 0 \\ 1.5753 & 0.7814 & 0.7814x & z_3\left(x\right) & z_4\left(x\right) & z_5\left(x\right) & z_6\left(x\right) \\ 0 & 0 & 0 & z_{3,xx}\left(x\right) & z_{4,xx}\left(x\right) & z_{5,xx}\left(x\right) & z_{6,xx}\left(x\right) \end{bmatrix}\mathbf{c} + \begin{Bmatrix} 0 \\ 25 \\ 19.535 \\ 0 \end{Bmatrix} x^2 + \begin{Bmatrix} 0 \\ 0 \\ 10192 \\ 39.070 \end{Bmatrix} \tag{A9}$$

$$Q\left(x\right) = \begin{bmatrix} 0 & 0 & 0.7814 & z_{3,x}\left(x\right) & z_{4,x}\left(x\right) & z_{5,x}\left(x\right) & z_{6,x}\left(x\right) \end{bmatrix}\mathbf{c} + 39.070x \tag{A10}$$

where

$$z_3\left(x\right) = e^{0.0307x}\cos\left(0.0237x\right) \tag{A11a}$$

$$z_4\left(x\right) = e^{0.0307x}\sin\left(0.0237x\right) \tag{A11b}$$

$$z_5\left(x\right) = e^{-0.0307x}\cos\left(0.0237x\right) \tag{A11c}$$

$$z_6\left(x\right) = e^{-0.0307x}\sin\left(0.0237x\right) \tag{A11d}$$

Using the above expressions of internal force fields, $\mathbf{X}$ and $\mathbf{U}$ in Eq. (74) can be obtained according to Eq. (75) as

$$\mathbf{X} = \begin{bmatrix} \mathbf{X}_{11} & \mathbf{X}_{12} \\ \mathbf{X}_{21} & \mathbf{X}_{22} \end{bmatrix}, \mathbf{U} = \begin{Bmatrix} \mathbf{0}_{6\times1} \\ \mathbf{U}_2 \end{Bmatrix} \tag{A12}$$

where

$$\mathbf{X}_{11} = \begin{bmatrix} 1 & 0 & 0 & 0 & 0 & 0 \\ 0 & 1 & 0 & 0 & 0 & 0 \\ 0 & 0 & 1 & 0 & 0 & 0 \\ 0 & 0 & 0 & 1 & 0 & 0 \\ 0 & 0 & 0 & 0 & 1 & 0 \\ 0 & 0 & 0 & -1 & 0 & 0 \end{bmatrix} \tag{A13a}$$

$$\mathbf{X}_{12} = \begin{bmatrix} 0 & 0 & 0 & 0 & 0 & 0 \\ 0 & 0 & 0 & 0 & 0 & 0 \\ 0 & 0 & 0 & 0 & 0 & 0 \\ 0 & 0 & 0 & 0 & 0 & 0 \\ 0 & 0 & 0 & 0 & 0 & 0 \\ 0 & 2.8552\text{E-}9 & 1.1200\text{E-}10 & 8.6771\text{E-}11 & -1.1200\text{E-}10 & 8.6771\text{E-}11 \end{bmatrix} \tag{A13b}$$

$$\mathbf{X}_{21}=\begin{bmatrix} 1 & 0 & 0 & 0 & 0 & 1.7644\text{E-6} \\ 0 & 1 & 2000 & 0 & 0 & 1.7026\text{E-5} \\ 0 & 0 & 1 & 0 & 0 & 1.7026\text{E-8} \\ 0 & 0 & 0 & 1 & 0 & -7.4777\text{E-21} \\ 0 & 0 & 0 & -2000 & 1 & 7.3184\text{E-19} \\ 0 & 0 & 0 & -1 & 0 & 7.4777\text{E-21} \end{bmatrix} \tag{A13c}$$

$$\mathbf{X}_{22}=\begin{bmatrix} -1.7026\text{E-8} & -1.7026\text{E-5} & 1.3094\text{E17} & -3.9518\text{E17} & 4.9938\text{E-11} & -9.8819\text{E-10} \\ -4.3610\text{E-7} & -2.9073\text{E-4} & -2.0164\text{E18} & 3.5510\text{E18} & 5.4120\text{E-7} & 5.1200\text{E-7} \\ -4.3610\text{E-10} & -4.3610\text{E-7} & -1.4613\text{E17} & 6.0963\text{E16} & 2.7130\text{E-10} & 2.6080\text{E-10} \\ -1.7205\text{E-22} & 6.7763\text{E-21} & -1.8718\text{E17} & 7.7737\text{E16} & 3.4781\text{E-10} & 3.3328\text{E-10} \\ -5.0822\text{E-19} & 5.7104\text{E-6} & 1.1552\text{E18} & -5.0979\text{E18} & -6.9746\text{E-7} & -6.5428\text{E-7} \\ 1.7205\text{E-22} & 2.8552\text{E-9} & 1.5647\text{E17} & -1.2883\text{E17} & -3.4781\text{E-10} & -3.3328\text{E-10} \end{bmatrix} \tag{A13d}$$

$$\mathbf{U}_2=\left\{-5.6647\text{E-1} \quad -7.0477 \quad -1.4316\text{E-2} \quad 2.8396\text{E-4} \quad 1.5642\text{E-3} \quad 1.5642\text{E-6}\right\}^{\mathrm{T}} \tag{A14}$$

By using Eqs. (77)-(79), the interpolation expressions for the four displacement components ($u(x)$, $w(x)$, $\theta(x)$ and $\phi(x)$) can be obtained. Then, the element stiffness matrix can be obtained according to Eq. (91) as

$$\mathbf{K}_e=\begin{bmatrix} \mathbf{K}_e^{11} & \mathbf{K}_e^{12} \\ \left(\mathbf{K}_e^{12}\right)^{\mathrm{T}} & \mathbf{K}_e^{22} \end{bmatrix} \tag{A15}$$

where

$$\mathbf{K}_e^{11}=\begin{bmatrix} 9.9727\text{E5} & -3.2503\text{E-2} & -3.9108\text{E7} & 1.5714\text{E8} & 2.4038\text{E6} & 1.6834\text{E8} \\ & 6.7493\text{E3} & 6.7493\text{E6} & -3.1066\text{E6} & 4.2433\text{E4} & 2.8132\text{E6} \\ & & 1.0576\text{E10} & -1.1285\text{E10} & -5.6128\text{E7} & -4.0908\text{E9} \\ & & & 2.6234\text{E11} & 4.5639\text{E9} & 3.3439\text{E11} \\ & & & & 8.8483\text{E7} & 6.1117\text{E9} \\ \text{sym.} & & & & & 4.4112\text{E11} \end{bmatrix} \tag{A16a}$$

$$\mathbf{K}_e^{12}=\begin{bmatrix} -9.9727\text{E5} & 3.2503\text{E-2} & 3.9108\text{E7} & -1.5715\text{E8} & 2.4038\text{E6} & -1.6838\text{E8} \\ 3.2503\text{E-2} & -6.7493\text{E3} & 6.7493\text{E6} & -3.1065\text{E6} & -4.2435\text{E4} & 2.8133\text{E6} \\ 3.9108\text{E7} & -6.7493\text{E6} & 2.9223\text{E9} & 5.0716\text{E9} & -1.4100\text{E8} & 9.7173\text{E9} \\ -1.5714\text{E8} & 3.1066\text{E6} & 5.0714\text{E9} & -1.3534\text{E11} & 2.4324\text{E9} & -1.7495\text{E11} \\ -2.4038\text{E6} & -4.2433\text{E4} & 1.4099\text{E8} & -2.4325\text{E9} & 4.2932\text{E7} & -3.0937\text{E9} \\ -1.6838\text{E6} & -2.8132\text{E6} & 9.7172\text{E9} & -1.7495\text{E11} & 3.0936\text{E9} & -2.2295\text{E11} \end{bmatrix} \tag{A16b}$$

$$\mathbf{K}_e^{22}=\begin{bmatrix} 9.9727\text{E5} & -3.2503\text{E-2} & -3.9108\text{E7} & 1.5715\text{E8} & -2.4038\text{E6} & 1.6838\text{E8} \\ & 6.7493\text{E3} & -6.7493\text{E6} & 3.1065\text{E6} & 4.2435\text{E4} & -2.8133\text{E6} \\ & & 1.0576\text{E10} & -1.1285\text{E10} & 5.6126\text{E7} & -4.0907\text{E9} \\ & & & 2.6235\text{E11} & -4.5639\text{E9} & 3.3439\text{E11} \\ & & & & 8.8483\text{E7} & -6.1117\text{E9} \\ \text{sym.} & & & & & 4.4112\text{E11} \end{bmatrix} \tag{A16c}$$

Meanwhile, the element equivalent nodal force vector can be obtained according to Eq. (92) as

$$\mathbf{F}_e=\begin{bmatrix} \left\{-3.4557\text{E4} \quad -5.0000\text{E4} \quad -1.5250\text{E7} \quad -3.8157\text{E7} \quad -9.4464\text{E5} \quad -6.6183\text{E7}\right\}^{\mathrm{T}} \\ \left\{3.4557\text{E4} \quad -5.0000\text{E4} \quad 1.5250\text{E7} \quad 3.8157\text{E7} \quad -9.4464\text{E5} \quad 6.6183\text{E7}\right\}^{\mathrm{T}} \end{bmatrix} \tag{A17}$$

For the SS beam, the boundary conditions can be set to

$$\begin{matrix} u^a & w^a & w^a_{,x} & \theta^a & \phi^a & \phi^a_{,x} & u^b & w^b & w^b_{,x} & \theta^b & \phi^b & \phi^b_{,x} \\ [\;1 & 1 & 0 & 0 & 1 & 0 & 0 & 1 & 0 & 0 & 1 & 0\;] \end{matrix} \tag{A18}$$

where '1' and '0' represent the constrained and free DoFs, respectively.

By introducing the boundary conditions, the solutions of nodal displacements can be obtained by using only one element as

$$\mathbf{\Phi} = \begin{bmatrix} \{0 \quad 0 \quad -3.7481\text{E-}3 \quad -1.4657\text{E-}4 \quad 0 \quad 7.0312\text{E-}6\}^{\mathrm{T}} \\ \{2.8480\text{E-}1 \quad 0 \quad 3.7481\text{E-}3 \quad 1.4657\text{E-}4 \quad 0 \quad -7.0312\text{E-}6\}^{\mathrm{T}} \end{bmatrix} \tag{A19}$$

Then, the mid-span vertical displacement can be obtained according to the interpolation expressions as

$$u_y\left(L/2,0\right) = w\left(L/2\right) + g\left(0\right)\phi\left(L/2\right) = \mathbf{N}^h_w\left(L/2\right)\mathbf{\Phi} + U^h_w\left(L/2\right) + \mathbf{N}^h_\phi\left(L/2\right)\mathbf{\Phi} + U^h_\phi\left(L/2\right) = 2.3268 \tag{A20}$$